\documentclass[onefignum,onetabnum]{siamart190516}
\usepackage[utf8]{inputenc}
\usepackage{latexsym,amsfonts,amsmath,amssymb,url}
\usepackage{mathtools}
\usepackage{xcolor,graphicx}
\usepackage{comment}

\newsiamthm{example}{Example}
\newsiamremark{assumption}{Assumption}
\newsiamremark{remark}{Remark}

\newcommand{\median}{\mathop{\rm median}\limits}
\newcommand{\ri}{\mathrm{i}}
\newcommand{\rd}{\, \mathrm{d}}
\newcommand{\bsh}{\boldsymbol{h}}
\newcommand{\bsk}{\boldsymbol{k}}

\newcommand{\bszero}{\boldsymbol{0}}
\newcommand{\bsalpha}{\boldsymbol{\alpha}}
\newcommand{\bsgamma}{\boldsymbol{\gamma}}
\newcommand{\bsomega}{\boldsymbol{\omega}}
\newcommand{\bstau}{\boldsymbol{\tau}}
\newcommand{\bsq}{\boldsymbol{q}}
\newcommand{\bsx}{\boldsymbol{x}}
\newcommand{\bsy}{\boldsymbol{y}}
\newcommand{\bsz}{\boldsymbol{z}}

\newcommand{\wal}{\mathrm{wal}}
\newcommand{\wor}{\mathrm{wor}}
\newcommand{\Kor}{\mathrm{kor}}
\newcommand{\Sob}{\mathrm{sob}}
\newcommand{\per}{\mathrm{per}}
\newcommand{\nonper}{\mathrm{nonper}}
\newcommand{\CC}{\mathbb{C}}

\newcommand{\FF}{\mathbb{F}}
\newcommand{\NN}{\mathbb{N}}
\newcommand{\PP}{\mathbb{P}}
\newcommand{\RR}{\mathbb{R}}
\newcommand{\UU}{\mathbb{U}}
\newcommand{\ZZ}{\mathbb{Z}}
\newcommand{\Scal}{\mathcal{S}}

\DeclareMathOperator{\tr}{tr}

\def\cF   {\mathcal{F}}

\def\cO   {\mathcal{O}}

\def\?{\discretionary{}{}{}}  

\allowdisplaybreaks

\newif\ifnotes\notestrue

\headers{Construction-free median QMC rules}{T. Goda and P. L'Ecuyer}

\title{Construction-free median quasi-Monte Carlo rules for function spaces with unspecified smoothness and general weights
\thanks{Submitted to the editors in January 2022.%
\funding{T.~Goda is supported by JSPS KAKENHI Grant Number 20K03744. P.~L'Ecuyer is supported by Discovery Grant number RGPIN-110050 from NSERC-Canada.}}}

\author{Takashi Goda\thanks{School of Engineering, University of Tokyo, 7-3-1 Hongo, Bunkyo-ku, Tokyo 113-8656, Japan (\email{goda@frcer.t.u-tokyo.ac.jp}).}
\and Pierre L'Ecuyer\thanks{DIRO, Universit\'{e} de Montr\'{e}al, C.P. 6128,
 Succ. Centre-Ville, Montr\'{e}al, Qc, Canada H3C 3J7, (\email{lecuyer@iro.umontreal.ca}).}}

\date{\today}

\begin{document}

\maketitle

\begin{abstract}
We study quasi-Monte Carlo (QMC) integration of smooth functions defined over the multi-dimensional unit cube. Inspired by a recent work of Pan and Owen, we study a new \emph{construction-free median QMC rule} which can exploit the smoothness and the weights of function spaces adaptively. For weighted Korobov spaces, we draw a sample of $r$ independent generating vectors of rank-1 lattice rules, compute the integral estimate for each, and approximate the true integral by the median of these $r$ estimates. For weighted Sobolev spaces, we use the same approach but with the rank-1 lattice rules replaced by high-order polynomial lattice rules. A major advantage over the existing approaches is that we do not need to construct good generating vectors by a computer search algorithm, while our median QMC rule achieves almost the optimal worst-case error rate for the respective function space with \emph{any smoothness and weights}, with a probability that converges to 1 exponentially fast as $r$ increases. Numerical experiments illustrate and support our theoretical findings.
\end{abstract}



\section{Introduction}
\label{sec:intro}

We consider numerical integration of functions defined over the $s$-dimensional unit cube $[0,1)^s$. For an integrable function $f: [0,1)^s\to \RR$, we denote the integral of $f$ by
\[ 
  I_s(f):=\int_{[0,1)^s}f(\bsx)\rd \bsx. 
\]
The \emph{quasi-Monte Carlo} (QMC) method approximates $I_s(f)$ by the equally-weighted average of function evaluations over a deterministic $N$-element point set $P_{N,s}\subset [0,1)^s$:
\[ 
  Q_{P_{N,s}}(f) = \frac{1}{N}\sum_{\bsx\in P_{N,s}}f(\bsx).
\]
The \emph{worst-case error} for a given normed function space $\cF$ and point set $P_{N,s}$ is 
\[
   e^{\wor}(Q_{P_{N,s}}; \cF) = \sup_{f\in \cF;\, \Vert f\Vert\le 1} \left| Q_{P_{N,s}}(f) - I_s(f) \right|
\]
where $\Vert f\Vert$ is the norm of $f$ in this space.  
The key to success of the QMC method lies in a proper choice of the point set depending on a target class of functions. 
One wants to construct point sets for which this worst-case error is small and converges at the 
fastest possible rate as a function of $N$, for the given space $\cF$.
In \emph{randomized QMC}, the point set $P_{N,s}$ is randomized in a way that $Q_{P_{N,s}}(f)$ becomes an 
unbiased estimator of $I_s(f)$ and one wishes to minimize its variance \cite{vLEC18a,vLEC02a,Lemieux09}.  
In this paper, we focus on deterministic QMC point sets only.

There are two main families of QMC point sets: digital nets and sequences \cite{DickPilli10,Niederreiter92} and lattice point sets \cite{Niederreiter92,SloanJoe94}. 
We refer the reader to \cite{DKS13,Lemieux09,LeobPilli14} for further introductory details. 
In this paper, we consider rank-1 lattice point sets for Korobov spaces of periodic functions, 
and high-order polynomial lattice point sets \cite{DP07,Niederreiter92b,Pillichshammer12}
(which are a special type of digital nets) for Sobolev spaces of non-periodic functions.
Each point set from these types is defined by an $s$-dimensional generating vector, with integer coordinates 
in the ordinary lattice case and with polynomial coordinates in the polynomial lattice case. 

In both cases, the weighted spaces of functions are defined by selecting a positive smoothness parameter $\alpha$
and a positive weight $\gamma_u$ for each subset of coordinates $u \subseteq\{1,\dots,s\}$, with $\gamma_{\emptyset}=1$.
The parameter $\alpha$ determines how smooth the admissible functions are required to be.
For the Korobov spaces, it tells the minimal rate at which the Fourier coefficients of $f$ are required to decay,
and when it is an integer it corresponds to the minimal number of square-integrable mixed partial derivatives 
of $f$ with respect to each coordinate; see \cite[Appendix~A]{NW_Tract1}.
For the Sobolev spaces, $\alpha$ is a positive integer which also imposes integrability conditions on the 
partial derivatives of $f$.
The weights $\gamma_u$ act as constant importance factors given to the subsets of coordinates \cite{DSWW06,SW98}.
A larger $\gamma_u$ means that the projection of $f$ over the subset of coordinates in $u$ can have a larger
variation in some sense, so that more importance should be given to the uniformity of the points over this projection.

It is known that the best possible QMC point sets cannot provide a better convergence rate than
$\cO(N^{-\alpha})$ for the worst-case error for these two function spaces.
On the other hand, there are effective search algorithms which, 
for a given $\alpha$ and a given selection of weights $\gamma_u$, 
can construct good rank-1 lattice or polynomial lattice point sets for the function spaces determined by these parameters,
and for which the worst-case error converges as $\cO(N^{-\alpha + \epsilon})$ for any $\epsilon > 0$ 
\cite{DickPilli10,Kuo03}.
Software that implement such algorithms is also freely available \cite{LatNetbuilder,Latticebuilder,iNUY20m}.
These algorithms typically use a greedy component-by-component (CBC) construction approach proposed originally by
\cite{Korobov59}, then re-introduced and popularized by \cite{SR02}. With the CBC approach, 
the generating vector is determined one coordinate at a time by optimizing a figure of merit that depends only on this new coordinate and the previous ones, and where the previous coordinates can no longer be changed.

In general, the number of weights $\gamma_u$ to specify is $2^s-1$. 
When $s$ is large, specifying all these weights individually
becomes impractical, so it is common practice to ``parameterize'' the weights by a smaller number of parameters, 
usually linear instead of exponential in $s$.  
The most popular forms of parameterizations are the product weights, the order-dependent weights, and their combination.
For the product weights, one specifies a weight $\gamma_j > 0$ for each coordinate $j=1,\dots,s$, 
and the $\gamma_u$'s are defined as $\gamma_u = \prod_{j\in u} \gamma_j$.  
For the order-dependent weights, $\gamma_u$ depends only on the cardinality of $u$:
$\gamma_u = \Gamma_{|u|}$, where $\Gamma_1, \Gamma_2, \dots, \Gamma_s$ are selected positive constants.
Their combination gives the product and order-dependent (POD) weights, for which 
$\gamma_u = \Gamma_{|u|} \prod_{j\in u} \gamma_j$ \cite{vNIC14a}.
The main reason why the most popular choices of weights have this form is that the existing search algorithms 
are truly efficient for large $s$ only when the weights have this specific POD form
\cite{DSWW06,KSS12,LatNetbuilder,Latticebuilder,vNIC14a,NC06}.
Then, by using a fast-CBC approach that speeds up the search by exploiting a fast Fourier transform 
\cite{Latticebuilder,NC06},
one can find a generating vector that gives a worst-case error of $\cO(N^{-\alpha + \epsilon})$ in 
$\cO(s N \log N)$ operations for ordinary lattices and $\cO(\alpha s N \log N)$ for polynomial lattices (with interlacing).

Although this form is convenient, the restriction to POD weights is limiting: for a given application,
the appropriate weights may be quite far from the POD form.  In this case, imposing POD implies that the 
point sets are constructed with the wrong weights. 
Moreover, even without constraints imposed on the form of weights, finding or approximating appropriate weights
and the appropriate $\alpha$ for a given application is generally very difficult \cite{vLEC12a}.  
When the points are constructed with the wrong weights, the QMC method can be quite ineffective in general.
These drawbacks have been addressed very partially in recent papers.
In \cite{EKNO21}, the authors introduces a construction algorithm that does not require the knowledge of $\alpha$.
It uses a CBC construction algorithm with a figure of merit that assumes $\alpha=1$, and for each coordinate it also 
constructs the generating vector one binary digit at a time.  
The method provides a convergence rate of $\cO(N^{-\alpha+\epsilon})$.
In \cite{DG21}, the authors study the stability of rank-1 lattice rules and polynomial lattice rules to a (limited) misspecification of $\alpha$ and the weights, for product and POD weights.  
They obtain worst-case error bounds for function spaces determined by parameters $\alpha'$ and $\bsgamma' = \{\gamma'_u\}$ 
when the rules are constructed using parameters $\alpha$ and $\bsgamma=\{\gamma_u\}$ instead, under certain conditions 
on those weights. These results are interesting but they do not completely eliminate the need to specify the weights.

The method studied in this paper requires no knowledge at all on $\alpha$ and the weights $\gamma_u$. 
No value needs to be specified for any of these parameters.
The algorithm is inspired by recent work from Pan and Owen \cite{PanOwen2021}, and works as follows.
For a fixed odd integer $r > 0$, we draw $r$ generating vectors independently and uniformly from the set 
of all admissible generating vectors.  For each of them, we compute the corresponding QMC approximation 
$Q_{P_{N,s}}(f)$, then we take the median $M(f)$ of these $r$ approximations as our final estimate of $I_s(f)$.
Since the method does not require the explicit construction of a good point set, we call it a
\emph{construction-free median QMC rule}.

Our main contribution is to prove that for $\cF$ representing either a weighted Korobov or weighted Sobolev space
determined by parameters $\alpha$ and $\bsgamma$,
the error $|M(f) - I_s(f)|$ obeys the following type of probabilistic bound:
For any $\epsilon > 0$ and $0 < \rho < 1$, there is a constant $c_1 = c_1(\alpha,\bsgamma,\epsilon) > 0$ 
(which depends on $\alpha$, the $\gamma_u$'s and $\epsilon$) such that
\[
	 \PP\left[ \sup_{f\in \cF;\, \Vert f\Vert\le 1} \left|M(f) - I_s(f) \right| \le \frac{c_1(\alpha,\bsgamma,\epsilon)}{(\rho N)^{\alpha-\epsilon}}\right] 
   \ge 1 - \rho^{(r+1)/2}/4.
\]
In other words, the worst-case error of the median estimator is bounded by a quantity that decreases 
almost at the best possible rate of $\cO(N^{-\alpha})$, with a probability that converges to 1 exponentially
fast as a function of $r$.
That is, we have a simple method that provides essentially the best possible convergence rate with very high probability,
without requiring any knowledge of $\alpha$ and the weights. 
The key reason why this is possible is that the vast majority of the choices of generating 
vectors turn out to be quite good and give a QMC approximation which is quite close to $I_s(f)$.  
Only a small minority give a large error.  For the vector giving the median value to be in that small
minority, there must be at least $(r+1)/2$ generating vectors in the sample of size $r$ that belong to this
small minority, and the probability that this happens decreases towards zero exponentially in $r$.  

The remainder is organized as follows.
In Section~\ref{sec:lattice}, we recall some basic facts on lattice rules for Korobov spaces,
and we prove our main result for the median estimator in this setting.   
In Section~\ref{sec:plr}, we do the same for high-order polynomial lattice rules in Sobolev spaces.
In Section~\ref{sec:experiments}, we report numerical experiments to support our theoretical findings.

\section{Lattice rules for Korobov spaces}
\label{sec:lattice}

\subsection{Definitions}
\label{sec:lattice-definitions}

Lattice point sets are well suited for performing numerical integration of smooth periodic functions. 
A rank-1 lattice point set is defined as follows:
\begin{definition}[rank-1 lattice point set]
\label{def:lattice}
Let $N\geq 2$ be the number of points and $\bsz = (z_1, \ldots, z_s) \in \{1,\ldots,N-1\}^s$. 
The rank-1 lattice point set defined by $N$ and the generating vector $\bsz$ is 
\[ 
  P_{N,s,\bsz} = \left\{ \left( \left\{ \frac{nz_1}{N}\right\}, \ldots, 
	 \left\{ \frac{nz_s}{N}\right\} \right)\in [0,1)^s \mid n=0,1,\ldots, N-1\right\}, 
\]
where $\{x\}:=x-\lfloor x\rfloor$ denotes the fractional part of a real $x \ge 0$. 
The QMC algorithm using $P_{N,s,\bsz}$ as a point set is called the rank-1 lattice rule with generating vector $\bsz$.
\end{definition}

Let $f: [0,1)^s\to \RR$ be periodic with an absolutely convergent Fourier series
\[ f(\bsx) = \sum_{\bsk\in \ZZ^s}\hat{f}(\bsk)\exp\left( 2\pi i \bsk\cdot \bsx\right),\]
where the dot product $\cdot$ denotes the usual inner product of two vectors on the Euclidean space $\RR^s$ and $\hat{f}(\bsk)$ denotes the $\bsk$-th Fourier coefficient of $f$:
\[ \hat{f}(\bsk) := \int_{[0,1)^s}f(\bsx)\exp\left( -2\pi i \bsk\cdot \bsx\right)\rd \bsx. \]
Note that $\hat{f}(\bszero)$ coincides with the integral $I_s(f)$. As a class of periodic functions, we consider the following \emph{weighted Korobov space}.

\begin{definition}[weighted Korobov space]\label{def:Kor_space}
Let $\alpha>1/2$ and $\bsgamma=\{\gamma_u\}_{u\subseteq \{1,\ldots,s\}}$ be a set of positive weights with $\gamma_{\emptyset}=1$. For a non-empty subset $u\subseteq \{1,\ldots,s\}$ and a vector $\bsk_u\in (\ZZ\setminus \{0\})^{|u|}$, we denote by $(\bsk_u,\bszero)$ the vector $\bsh\in \ZZ^s$ such that $h_j=k_j$ if $j\in u$ and $h_j=0$ otherwise, and define
\begin{align*} r_{\alpha,\bsgamma}(\bsk_u,\bszero):=\gamma_u\prod_{j\in u}\frac{1}{|k_j|^{\alpha}},
\end{align*}
and set $r_{\alpha,\bsgamma}(\bszero)=1$. The weighted Korobov space, denoted by $\cF^{\Kor}_{s,\alpha,\bsgamma}$, is a reproducing kernel Hilbert space with reproducing kernel
\[ 
 K^{\Kor}_{s,\alpha,\bsgamma}(\bsx,\bsy) = \sum_{\bsk\in \ZZ^s}(r_{\alpha,\bsgamma}(\bsk))^2\exp\left( 2\pi i \bsk\cdot (\bsx-\bsy)\right), 
\]
and inner product
\[ 
  \langle f,g\rangle^{\Kor}_{s,\alpha,\bsgamma} = \sum_{\bsk\in \ZZ^s}\frac{\hat{f}(\bsk)\overline{\hat{g}(\bsk)}}{(r_{\alpha,\bsgamma}(\bsk))^2}.
\] 
We denote the induced norm by $\|f\|^{\Kor}_{s,\alpha,\bsgamma}:=\sqrt{\langle f,f\rangle^{\Kor}_{s,\alpha,\bsgamma}}$.
\end{definition}

One wishes to have a good generating vector $\bsz$ such that the worst-case error of the corresponding lattice rule for $\cF^{\Kor}_{s,\alpha,\bsgamma}$, defined by
\[
  e^{\wor}(Q_{P_{N,s,\bsz}};\cF^{\Kor}_{s,\alpha,\bsgamma}) 
	:=  \sup_{\substack{f\in \cF^{\Kor}_{s,\alpha,\bsgamma}\\ 
	 \|f\|^{\Kor}_{s,\alpha,\bsgamma}\leq 1}}|Q_{P_{N,s,\bsz}}(f)-I_s(f)|, 
\]
is small.  No good explicit construction scheme for such a $\bsz$ is known for $s\geq 3$, 
so that we usually resort to a computer search algorithm as mentioned earlier.
By restricting each $z_j$ to be in the set
\[ 
  \UU_N:=\{1\leq z\leq N-1\mid \gcd(z,N)=1\},
\]
we ensure that each projection of $P_{N,s,\bsz}$ on a single coordinate contains 
the $N$ distinct values $\{0, 1/N, \dots, (N-1)/N\}$ (no superposed points).
The CBC construction algorithm for a good generating vector $\bsz$ starts with $z_1=1$, 
then for $j=2,\dots,s$ it searches for the best component $z_j$ from the set $\UU_N$ 
while keeping the earlier components $z_1,\ldots,z_{j-1}$ unchanged.

For our median QMC rank-1 lattice rule for weighted Korobov spaces, we select an independent random sample
$\bsz_1,\dots,\bsz_r$ from the set $\UU_N^s$, 
and we approximate $I_s(f)$ by the median
\[ 
	 M_{N,s,r}(f) := \median\left(Q_{P_{N,s,\bsz_1}}(f),\ldots,Q_{P_{N,s,\bsz_r}}(f)\right).
\]
Note that, for given $\bsz_1,\dots,\bsz_r$, the index $\ell$ for which $\bsz_\ell$ gives the median $M_{N,s,r}(f)$ generally depends on $f$.
The worst-case error in this case is the random variable
\begin{equation}
\label{eq:ewor}
  e^{\wor}(M_{N,s,r};\cF^{\Kor}_{s,\alpha,\bsgamma}) 
	:=  \sup_{\substack{f\in \cF^{\Kor}_{s,\alpha,\bsgamma}\\ 
	 \|f\|^{\Kor}_{s,\alpha,\bsgamma}\leq 1}} |M_{N,s,r}(f) -I_s(f)|. 
\end{equation}
In this random expression, we assume that $\bsz_1,\dots,\bsz_r$ are first picked randomly,
then $f$ is taken as the worst-case function for the median, for these given $\bsz_1,\dots,\bsz_r$.

\subsection{Our main results on lattice rules for Korobov spaces}

To prove our main result, we need a few more definitions.

\begin{definition}[dual lattice]\label{def:dual}
For $N \geq 2$ and $\bsz \in \UU_N^s$, the set
\[ 
  P^{\perp}_{N,s,\bsz} := \left\{ \bsk\in \ZZ^s \mid \bsk\cdot\bsz\equiv 0 \pmod N\right\} 
\]
is called the dual lattice of the rank-1 lattice point set $P_{N,s,\bsz}$.
\end{definition}
\noindent The following character property of the rank-1 lattice rule is well-known, 
see for instance \cite[Lemmas~4.2 and 4.3]{DHP15}.
\begin{lemma}[character property]\label{lem:character}
For $N \geq 2$, $\bsz \in \UU_N^s$ and $\bsk\in \ZZ^s$, we have
\[ 
  \frac{1}{N}\sum_{\bsx\in P_{N,s,\bsz}}\exp\left( 2\pi i \bsk\cdot \bsx\right)=\begin{cases} 1 & \text{if $\bsk\in P^{\perp}_{N,s,\bsz}$,} \\ 
	0 & \text{otherwise.} \end{cases}
\]
\end{lemma}

As our first main result, we prove a probabilistic upper bound on the worst-case error of our 
median rank-1 lattice rule for weighted Korobov spaces.

\begin{theorem}\label{thm:main1}
Let $N\geq 2$ be an integer, $r > 0$ be an odd integer, and $\bsz_1,\ldots,\bsz_r$ be chosen independently and randomly from the set $\UU_N^s$ (with replacement). 
Then, for any $\alpha>1/2$ and $\bsgamma$, the worst-case error of the median rule obeys the following bound:
\[
 e^{\wor}(M_{N,s,r}; \cF^{\Kor}_{s,\alpha,\bsgamma})
 \leq \inf_{1/(2\alpha) < \lambda<1}\left(\frac{1}{\eta\varphi(N)}
  \sum_{\emptyset \neq u\subseteq \{1,\ldots,s\}}\gamma_u^{2\lambda}(2\zeta(2\alpha\lambda))^{|u|}\right)^{1/(2\lambda)}
\]
with a probability of at least
\[ 1-\binom{r}{(r+1)/2}\eta^{(r+1)/2},\]
for any $0<\eta<1$, where $\varphi$ and $\zeta$ denote the Euler totient function and the Riemann zeta function, respectively.
\end{theorem}

We note that the result for $r=1$, i.e., the case without taking the median, can be found, for instance, in \cite[Theorem~2]{DSWW06}, and has been used together with a random choice of $N$ in \cite{KKNU19} to prove an improved rate of convergence of the randomized error.

The following inequality on medians is a key ingredient in the proof of the theorem. 
Although it can be regarded as a special case of Jensen's inequality on medians proven in \cite{Merkle05}, 
we give a short direct proof to make the paper more self-contained.

\begin{lemma}\label{lem:median_jensen}
For any odd integer $r$ and real numbers $a_1,\ldots,a_r$, it holds that
\[ \left|\median\left(a_1,\ldots,a_r\right)\right| \leq \median\left(|a_1|,\ldots,|a_r|\right).\]
\end{lemma}
\begin{proof}
Because $r$ is odd, the median is unique. Let $\median(a_1,\ldots,a_r)=a_m$ for some $m\in \{1,\ldots,r\}$.
If $a_m \ge 0$, then $|a_m| = a_m \leq \median (|a_1|,\? \ldots,\? |a_r|)$.
If $a_m < 0$, there are at least $(r-1)/2$ other $a_\ell$'s for which $a_\ell \le a_m < 0$,
so $|a_\ell| \ge |a_m| > 0$.  
Then, $|\median(a_1,\ldots,a_r)| = |a_m| \le \median\left(|a_1|,\ldots,|a_r|\right)$.
\end{proof}

\noindent The same inequality holds for even $r$ by defining the median to be the arithmetic mean of the two middle values. However, in this paper, we focus on the case where $r$ is odd for the sake of simplicity.

We now prove our main result.
\begin{proof}[Proof of Theorem~\ref{thm:main1}]
Since any $f\in \cF^{\Kor}_{s,\alpha,\bsgamma}$ has an absolutely convergent Fourier series, by applying Lemma~\ref{lem:character}, Lemma~\ref{lem:median_jensen} and the Cauchy–Schwarz inequality, it holds for given $\bsz_1,\ldots,\bsz_r$ that
\begin{align}  
    & \qquad e^{\wor}(M_{N,s,r}; \cF^{\Kor}_{s,\alpha,\bsgamma}) \notag\\
	& = \sup_{\substack{f\in \cF^{\Kor}_{s,\alpha,\bsgamma} \notag \\ \|f\|^{\Kor}_{s,\alpha,\bsgamma}\leq 1}}\left| \median_{1\leq \ell\leq r}\frac{1}{N}\sum_{\bsx\in P_{N,s,\bsz_{\ell}}}f(\bsx)-I(f)\right| \notag \\
    & = \sup_{\substack{f\in \cF^{\Kor}_{s,\alpha,\bsgamma}\\ \|f\|^{\Kor}_{s,\alpha,\bsgamma}\leq 1}}\left| \median_{1\leq \ell\leq r}\frac{1}{N}\sum_{\bsx\in P_{N,s,\bsz_{\ell}}}\sum_{\bsk\in \ZZ^s}\hat{f}(\bsk)\exp\left( 2\pi i \bsk\cdot \bsx\right)-\hat{f}(\bszero)\right| \notag \\
    & = \sup_{\substack{f\in \cF^{\Kor}_{s,\alpha,\bsgamma}\\ \|f\|^{\Kor}_{s,\alpha,\bsgamma}\leq 1}}\left| \median_{1\leq \ell\leq r}\sum_{\bsk\in P^{\perp}_{N,s,\bsz_{\ell}}\setminus \{\bszero\}}\hat{f}(\bsk)\right| \notag \\
    & \leq \sup_{\substack{f\in \cF^{\Kor}_{s,\alpha,\bsgamma}\\ \|f\|^{\Kor}_{s,\alpha,\bsgamma}\leq 1}} \median_{1\leq \ell\leq r}\sum_{\bsk\in P^{\perp}_{N,s,\bsz_{\ell}}\setminus \{\bszero\}}|\hat{f}(\bsk)| \notag \\
    & = \sup_{\substack{f\in \cF^{\Kor}_{s,\alpha,\bsgamma}\\ \|f\|^{\Kor}_{s,\alpha,\bsgamma}\leq 1}} \median_{1\leq \ell\leq r}\sum_{\bsk\in P^{\perp}_{N,s,\bsz_{\ell}}\setminus \{\bszero\}}\frac{|\hat{f}(\bsk)|}{r_{\alpha,\bsgamma}(\bsk)}r_{\alpha,\bsgamma}(\bsk) \notag \\
    & \leq  \sup_{\substack{f\in \cF^{\Kor}_{s,\alpha,\bsgamma}\\ \|f\|^{\Kor}_{s,\alpha,\bsgamma}\leq 1}} \median_{1\leq \ell\leq r}\left(\sum_{\bsk\in P^{\perp}_{N,s,\bsz_{\ell}}\setminus \{\bszero\}}\frac{|\hat{f}(\bsk)|^2}{(r_{\alpha,\bsgamma}(\bsk))^2}\right)^{1/2} \notag \\ 
    & \qquad \times \left(\sum_{\bsk\in P^{\perp}_{N,s,\bsz_{\ell}}\setminus \{\bszero\}}(r_{\alpha,\bsgamma}(\bsk))^2\right)^{1/2} \notag \\
    & \leq \sup_{\substack{f\in \cF^{\Kor}_{s,\alpha,\bsgamma}\\ \|f\|^{\Kor}_{s,\alpha,\bsgamma}\leq 1}}\left(\sum_{\bsk\in \ZZ^s\setminus \{\bszero\}}\frac{|\hat{f}(\bsk)|^2}{(r_{\alpha,\bsgamma}(\bsk))^2}\right)^{1/2} 
	\median_{1\leq \ell\leq r} \Scal_{\alpha,\bsgamma}(\bsz_\ell) \notag \\
    & \leq \median_{1\leq \ell\leq r} \Scal_{\alpha,\bsgamma}(\bsz_\ell),\label{eq:proof_main_1}
\end{align}
in which
\[ 
  \Scal_{\alpha,\bsgamma}(\bsz) 
	:= \left(\sum_{\bsk\in P_{N,s,\bsz}^{\perp}\setminus \{\bszero\}}(r_{\alpha,\bsgamma}(\bsk))^2\right)^{1/2}.
\]

For $1/(2\alpha)<\lambda\leq 1$, by using the subadditivity
\begin{align}\label{eq:jensen} 
 \left( \sum_i a_i\right)^{\lambda}\leq \sum_i a_i^{\lambda}, 
\end{align}
which holds for non-negative reals $a_1,a_2,\ldots>0$, see \cite[Theorem~2.2]{DHP15}, and noting that the cardinality of $\UU_N$ is equal to $\varphi(N)$, we have
\begin{align*}
    & \ \ \ \ \ \frac{1}{(\varphi(N))^s}\sum_{\bsz\in \UU_N^s}(\Scal_{\alpha,\bsgamma}(\bsz))^{2\lambda} \\
	& = \ \frac{1}{(\varphi(N))^s}\sum_{\bsz\in \UU_N^s}\left(\sum_{\bsk\in P_{N,s,\bsz}^{\perp}\setminus \{\bszero\}}(r_{\alpha,\bsgamma}(\bsk))^2\right)^{\lambda} \\
    & \leq \ \frac{1}{(\varphi(N))^s}\sum_{\bsz\in \UU_N^s}\sum_{\bsk\in P_{N,s,\bsz}^{\perp}\setminus \{\bszero\}}(r_{\alpha,\bsgamma}(\bsk))^{2\lambda}\\
    & = \ \sum_{\bsk\in \ZZ^s\setminus \{\bszero\}}(r_{\alpha,\bsgamma}(\bsk))^{2\lambda}\frac{1}{(\varphi(N))^s}\sum_{\bsz\in \UU_N^s}\frac{1}{N}\sum_{\bsx\in P_{N,s,\bsz}}\exp(2\pi i\bsk\cdot \bsx)\\
    & = \ \sum_{\bsk\in \ZZ^s\setminus \{\bszero\}}(r_{\alpha,\bsgamma}(\bsk))^{2\lambda}\frac{1}{(\varphi(N))^s}\sum_{\bsz\in \UU_N^s}\frac{1}{N}\sum_{n=0}^{N-1}\exp(2\pi i n\bsk\cdot \bsz /N)\\
    & = \ \sum_{\emptyset\neq u\subseteq \{1,\ldots,s\}}\gamma_u^{2\lambda}\sum_{\bsk_u\in (\ZZ\setminus \{0\})^{|u|}}\frac{1}{(\varphi(N))^s} \sum_{\bsz\in \UU_N^s}\frac{1}{N}\sum_{n=0}^{N-1}\prod_{j\in u}\frac{\exp(2\pi i n k_j z_j /N)}{|k_j|^{2\alpha \lambda}}\\
    & = \ \frac{1}{N}\sum_{n=0}^{N-1}\sum_{\emptyset\neq u\subseteq \{1,\ldots,s\}}\frac{\gamma_u^{2\lambda}}{(\varphi(N))^{|u|}}\sum_{\bsz_u\in \UU_N^{|u|}}\sum_{\bsk_u\in (\ZZ\setminus \{0\})^{|u|}}\prod_{j\in u}\frac{\exp(2\pi i n k_jz_j/N)}{|k_j|^{2\alpha \lambda}}\\
    & = \ \frac{1}{N}\sum_{n=0}^{N-1}\sum_{\emptyset\neq u\subseteq \{1,\ldots,s\}}\gamma_u^{2\lambda}(T_{2\alpha \lambda}(n,N))^{|u|},
\end{align*}
where we write
\[ T_{2\alpha \lambda}(n,N):= \frac{1}{\varphi(N)}\sum_{z\in \UU_N}\sum_{k\in \ZZ\setminus \{0\}}\frac{\exp(2\pi i n kz/N)}{|k|^{2\alpha \lambda}}.\]

Since it follows from \cite[Lemmas~2.1 \& 2.2]{KJ02} that, for any positive integer $d$
\[ \frac{1}{N}\sum_{n=0}^{N-1}(T_{2\alpha \lambda}(n,N))^d\leq \frac{(2\zeta(2\alpha\lambda))^d}{\varphi(N)},\]
we obtain
\begin{align*}
    \frac{1}{(\varphi(N))^s}\sum_{\bsz\in \UU_N^s}(\Scal_{\alpha,\bsgamma}(\bsz))^{2\lambda} & \leq \sum_{\emptyset\neq u\subseteq \{1,\ldots,s\}}\gamma_u^{2\lambda}\frac{1}{N}\sum_{n=0}^{N-1}(T_{2\alpha \lambda}(n,N))^{|u|} \\
    & \leq \frac{1}{\varphi(N)}\sum_{\emptyset\neq u\subseteq \{1,\ldots,s\}}\gamma_u^{2\lambda}(2\zeta(2\alpha\lambda))^{|u|}.
\end{align*}
This gives an upper bound on the average of $(\Scal_{\alpha,\bsgamma}(\bsz))^{2\lambda}$ over all of the admissible 
$\bsz\in \UU_N^s$, which holds for any $1/(2\alpha)<\lambda\leq 1$. 

Then, Markov's inequality ensures that for any $0<\eta<1$, the probability of having
\[ 
 \Scal_{\alpha,\bsgamma}(\bsz) > \inf_{1/(2\alpha)<\lambda<1}\left(\frac{1}{\eta\varphi(N)}\sum_{\emptyset \neq u\subseteq \{1,\ldots,s\}}\gamma_u^{2\lambda}(2\zeta(2\alpha\lambda))^{|u|}\right)^{1/(2\lambda)}  =: B(\alpha,\bsgamma)
\]
is at most $\eta$ for a random choice of $\bsz\in \UU_N^s$. 
For the median estimator $M_{N,s,r}$ to be larger than this bound $B(\alpha,\bsgamma)$,
we must have $\Scal_{\alpha,\bsgamma}(\bsz_\ell) > B(\alpha,\bsgamma)$ for at least 
$(r+1)/2$ vectors among $\bsz_1,\ldots,\bsz_r$.
Taking the union bound on possible sets of $(r+1)/2$ vectors with $\Scal_{\alpha,\bsgamma}(\bsz_\ell) > B(\alpha,\bsgamma)$, the probability that this happens is bounded above by
\[ 
  \binom{r}{(r+1)/2}\eta^{(r+1)/2}.
\]
Combining this with the bound shown in \eqref{eq:proof_main_1} completes the proof.
\end{proof}

\begin{remark}\label{rem:probability}
One can easily prove by induction on $k$ that $\binom{2k-1}{k} < 4^{k-1}$ for $k \ge 2$.
Indeed, this is true for $k=2$, and for $k\ge 2$, one has
\[
  \binom{2k+1}{k+1} = \frac{2(2k+1)}{(k+1)} \binom{2k-1}{k} < 4\binom{2k-1}{k}<4^k.
\]
Then, for any odd $r\ge 3$, we have
\begin{equation}
\label{eq:binom-bound}
  \binom{r}{(r+1)/2} \eta^{(r+1)/2} < (4\eta)^{(r+1)/2}/4.
\end{equation}
Thus, the probability given in Theorem~\ref{thm:main1} must be larger than $1-(4\eta)^{(r+1)/2}/4$, 
which converges to 1 exponentially fast as a function of $r$ for $0<\eta<1/4$.
\end{remark}

By taking $1/(2\lambda) = \alpha-\epsilon$ and $\rho = 4\eta$ and using the previous remark,
we obtain the following corollary as a simplified version of Theorem~\ref{thm:main1}.

\begin{corollary}
For any odd $r\ge 3$, $\epsilon > 0$, and $0 < \rho < 1$, there is a constant $c_1 = c_1(\alpha,\bsgamma,\epsilon) > 0$ 
(which depends on $\alpha$, the $\gamma_u$'s and $\epsilon$) such that
\[
  \PP\left[e^{\wor}(M_{N,s,r}; \cF^{\Kor}_{s,\alpha,\bsgamma})
	         \le \frac{c_1(\alpha,\bsgamma,\epsilon)}{(\rho N)^{\alpha-\epsilon}}\right] 
  \ge 1 - \rho^{(r+1)/2}/4.
\]
\end{corollary}

\begin{proof}
For prime $N$, we have $\varphi(N) = N-1$ and the corollary follows 
from Theorem~\ref{thm:main1} and the bound (\ref{eq:binom-bound}) in Remark~\ref{rem:probability}.
For a general $N$, we know from  \cite[Theorem~15]{RS62} that 
\[ 
  \frac{1}{\varphi(N)}\leq \frac{1}{N}\left[ e^C\log\log N+\frac{2.50637}{\log\log N}\right]
\]
for any $N\geq 3$, where $C=0.57721\ldots$ is the Euler's constant.
From this, the probabilistic bound follows.
\end{proof}

\begin{remark}\label{rem:tractability}
Under some additional conditions on the weights $\bsgamma$, this bound depends only polynomially on the dimension $s$,
and can even be independent on the dimension $s$. We refer to \cite[Theorem~3]{DSWW06} for the case of general weights and \cite[Theorem~4]{DSWW06} for the case of product weights.
\end{remark}

\begin{remark}
As mentioned in Section~\ref{sec:intro}, our median rank-1 lattice rule is motivated by the observation that most of the possible generating vectors $\bsz \in \UU_N^s$ are a good choice, but the remaining ones are bad. To show this, let us first point out that $\Scal_{\alpha,\bsgamma}(\bsz)$ coincides with the worst-case error of the rank-1 lattice rule with the given $\bsz$, see \cite[Theorem~5.12]{DKS13}. As already seen in the proof of Theorem~\ref{thm:main1}, a proportion of the generating vectors $\bsz$ which satisfy the bound of order $N^{-\alpha+\epsilon}$, i.e.,
\[ \Scal_{\alpha,\bsgamma}(\bsz)\leq \inf_{1/(2\alpha)<\lambda<1}\left(\frac{1}{\eta\varphi(N)}\sum_{\emptyset \neq u\subseteq \{1,\ldots,s\}}\gamma_u^{2\lambda}(2\zeta(2\alpha\lambda))^{|u|}\right)^{1/(2\lambda)}, \]
is greater than or equal to $1-\eta$, for any $0< \eta <1$. On the other hand, the averaging argument in the proof of Theorem~\ref{thm:main1} with $\lambda=1$ gives
\[ 
 \frac{1}{(\varphi(N))^s}\sum_{\bsz\in \UU_N^s}(\Scal_{\alpha,\bsgamma}(\bsz))^2
  \leq \frac{1}{\varphi(N)}\sum_{\emptyset\neq u\subseteq \{1,\ldots,s\}}\gamma_u^2(2\zeta(2\alpha))^{|u|}.
\]
This implies that, for each $N$, there exists a small, distinct set of ``bad'' generating vectors whose $\Scal_{\alpha,\bsgamma}(\bsz)$ values are quite large so that the average of the \emph{squared} worst-case error over all the possible generating vectors is merely of order $N^{-1}$.
Such bad vectors may include those with all components being the same.
\end{remark}

\begin{figure}[hbt]
    \centering
    \includegraphics[width=0.45\linewidth]{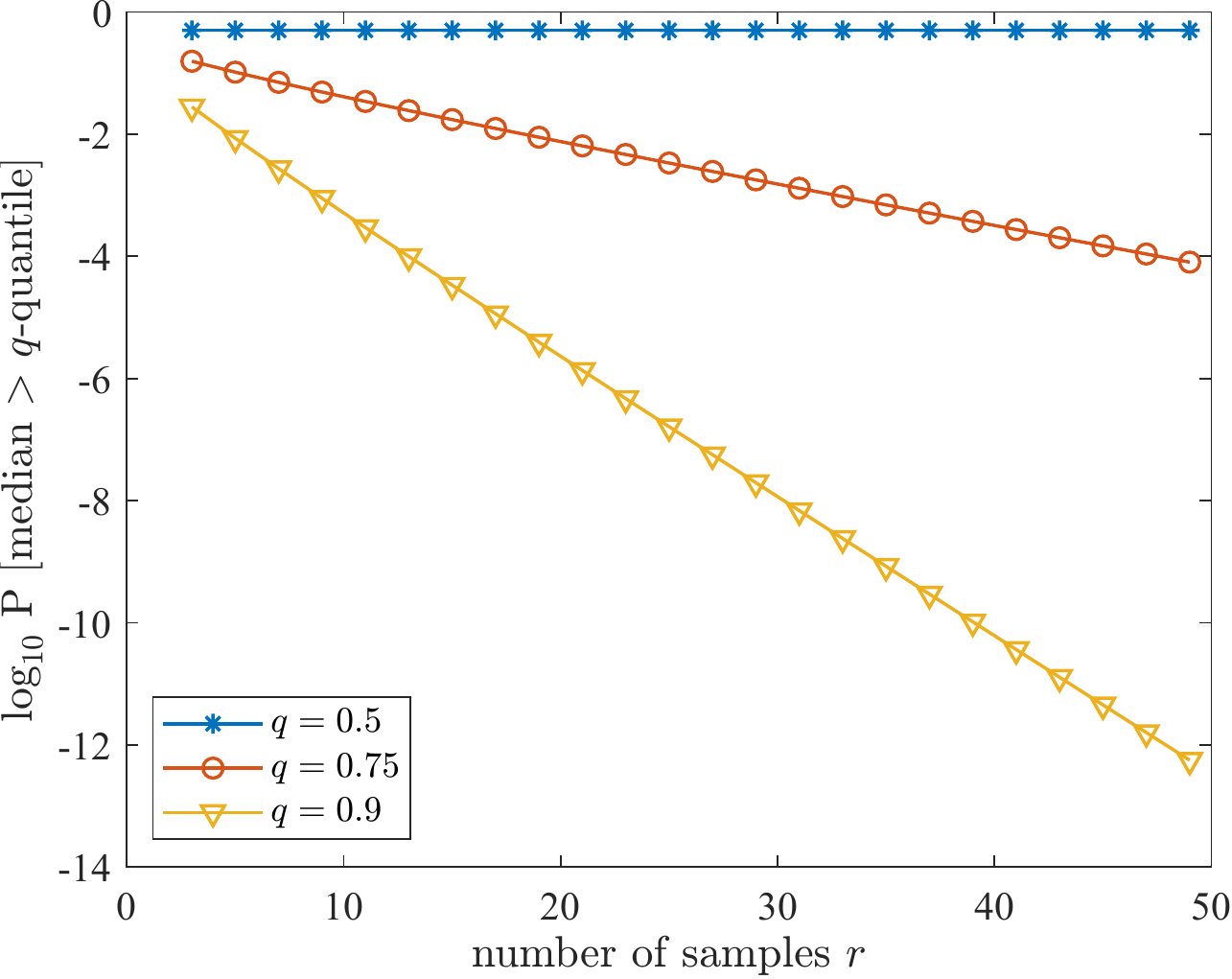}\\
    \caption{Probability $p_{+}(r,q)$ (on a $\log_{10}$ scale) as a function of $r$ for $q=0.5$ (blue), $q=0.75$ (red) and $q=0.9$ (yellow).}
\label{fig:quantile_prob}
\end{figure}

\begin{remark}
\label{rem:binomial}
Theorem~\ref{thm:main1} (or its corollary) gives a probabilistic error bound on the worst-case error (\ref{eq:ewor}),
together with a lower bound on the probability that this error bound holds.
The exact value of this probability depends on the probability distribution of $Q_{P_{N,s,\bsz}}(f)$
when $\bsz$ is drawn uniformly from $\UU_N^s$, 
and on the choices of $r$ and of the other parameters in the error bound. 
To get some insight on how it behaves, we will simplify the setting slightly and look at a one-side
error bound for a fixed $f$: we want to estimate the probability that the median $M_{N,s,r}(f)$ 
does not exceed some arbitrary constant $y$ larger than the mean $I_s(f)$.  
Suppose that this $y$ is the $q$-quantile $y_q$ of the distribution of $Q_{P_{N,s,\bsz}}(f)$ 
for some $q\in (3/4,\,1)$, i.e., $q = \PP[Q_{P_{N,s,\bsz}}(f) \le y_q]$.
Then the median $M_{N,s,r}(f)$ is larger than $y_q$ if and only if 
at least $(r+1)/2$ values are larger than $y_q$, and the probability that this happens is
\begin{align}
\label{eq:quantile_prob_right}
  p_{+}(r,q) = \sum_{i=(r+1)/2}^{r} \binom{r}{i} (1-q)^i q^{r-i}.  
\end{align}
Figure~\ref{fig:quantile_prob} plots $\log_{10} p_{+}(r,q)$ as a function of 
$r\in \{3,5,\ldots,49\}$ for $q=0.5$, $q=0.75$ and $q=0.9$, respectively. 
We see that $p_{+}(r,0.5)=0.5$ for any $r$ and that $p_{+}(r,q) \approx 10^{-\gamma r}$ where $\gamma \approx 0.071$ for $q=0.75$ and 
$\gamma \approx 0.231$ for $q=0.9.$
These plots provide some insight on the choice of $r$.
In particular, for fixed $q>0.5$, doubling $r$ squares the probability $p_{+}(r,q)$.
Suppose for example that we want $p_{+}(r,q) \le 10^{-4}$, to have a reasonable assurance 
that $M_{N,s,r}(f) \le y_q$.  
The plot shows that the minimal value of $r$ for this is about
$r = 13$ for $q = 0.9$, and about $r=49$ for $q= 0.75$.
For a given $f \in \cF^{\Kor}_{s,\alpha,\bsgamma}$ and fixed $N$, a larger $q$ means a larger $y_q$, 
but for a fixed $q$ we can reduce $y_q$ and bring it close to $I_s(f)$ by increasing $N$. 
From Theorem~\ref{thm:main1} with $r=1$ and $\eta = 1-q < 1/4$ (or its corollary with $\rho = 4(1-q) < 1$), 
we have that $|y_q - I_s(f)|$ is $\cO(N^{-\alpha + \epsilon})$.
In summary, for a fixed $q > 3/4$, we can decrease the error bound by increasing $N$
and increase the probability that the bound holds by increasing $r$.  
We can also increase both $q$ and $N$ in a way such that $y_q$ remains about the same; 
then the same $p_{+}(r,q)$ can be obtained with a reduced $r$.
What we just said is for the upper bound $M_{N,s,r}(f) \le y_q$, but essentially the same discussion 
can be made concerning the assurance that $M_{N,s,r}(f) > y_{1-q}$.
In applications, the values of $y_q$ and $y_{1-q}$ are unknown, but our reasoning 
suggests that a moderate value of $r$, say no more than 25, should be sufficient in practice,
together with a large $N$ (as large as the computing budget allows). 
The results of our numerical experiments support this.
\end{remark}

\begin{remark}
It is known that rank-1 lattice rules also work for non-periodic functions by applying the tent transformation
\[ \pi(x)=1-|2x-1|\]
component-wise to every point in the set $P_{N,s,\bsz}$ \cite{Hickernell02,DNP14,CKNS16,GSY19}. The same probabilistic upper bound, shown in Theorem~\ref{thm:main1}, holds for the worst-case error of the median rule built up of the tent-transformed rank-1 lattice rules in the so-called \emph{weighted half-period cosine spaces} with any parameter $\alpha$ and weights $\bsgamma$. As shown in \cite[Lemma~1]{DNP14}, the half-period cosine space coincides with an unanchored Sobolev space with smoothness 1 when $\alpha=1$.
\end{remark}

\section{High-order polynomial lattice rules for Sobolev spaces}
\label{sec:plr}

We now consider high-order polynomial lattice point sets as defined in \cite{DP07}. 
These point sets are well suited for performing numerical integration of smooth non-periodic functions. In what follows let $b$ be a prime, and $\FF_b$ be the finite field of order $b$, which we identify with the set $\{0,1,\ldots,b-1\}$. Let $\NN$ be the set of positive integers and $\NN_0:=\NN\cup \{0\}$. For $k\in \NN_0$ having the $b$-adic finite expansion $k=\kappa_0+\kappa_1 b+\cdots$, we write $k(x)=\kappa_0+\kappa_1 x+\cdots \in \FF_b[x]$. 
With these ingredients, we have the following definition from \cite{DP07}.

\begin{definition}[high-order polynomial lattice point sets]\label{def:HO_poly_lattice}
Let $m, n\in \NN$ with $m\le n$, $p\in \FF_b[x]$ with $\deg(p)=n$ and $\bsq = (q_1, \ldots, q_s) \in (\FF_b[x])^s$ with $\deg(q_j)<n$. The high-order polynomial lattice point set defined by $m,n,p$ and $\bsq$ consists of $N=b^m$ points and is given by
\[\small 
  P_{m,n,s,p,\bsq} = \left\{ \left( \nu_n\left( \frac{h(x)q_1(x)}{p(x)}\right), \ldots, \nu_n\left( \frac{h(x)q_s(x)}{p(x)}\right) \right)\in [0,1)^s \mid h=0,\ldots, b^m-1\right\}, 
\]
where $\nu_n: \FF_b((x^{-1}))\to [0,1)$ is defined by
\[ 
 \nu_n\left( \sum_{i=w}^{\infty}\frac{a_i}{x^i}\right) := \sum_{i=\max(1,w)}^{n}\frac{a_i}{b^i}.
\]
The QMC algorithm using $P_{m,n,s,p,\bsq}$ as a point set is called the high-order polynomial lattice rule with modulus $p$ and generating vector $\bsq$.   The \emph{order} of this rule is defined as $\lfloor n/m\rfloor$.
\end{definition}

Typically, $n$ will be a multiple of $m$.  When $n=m$, this gives the digital net construction 
introduced in \cite{Niederreiter92b} and called polynomial lattice rule in \cite{DickPilli10}.
Note that \cite{vLEC02a,vLEM03a} introduced the term ``polynomial lattice rule'' with a slightly different definition,
in which the coordinates of the points have an infinite periodic expansion and the modulus $p$ has degree $m$.
The construction in Definition~\ref{def:HO_poly_lattice} essentially builds a polynomial lattice point set 
with $b^n$ points and uses only the first $b^m$ points.

Instead of the weighted Korobov space $\cF^{\Kor}_{s,\alpha,\bsgamma}$, we consider the following Sobolev-type Banach space as our target space for high-order polynomial lattice rules. 

\begin{definition}[weighted Sobolev space]\label{def:Sob_space}
Let $\alpha\in \NN$, $\alpha\geq 2$, $1\leq q\leq \infty$ and let $\bsgamma=\{\gamma_u\}_{u\subseteq \{1,\ldots,s\}}$ be a set of positive weights with $\gamma_{\emptyset}=1$. The weighted Sobolev space, denoted by $\cF^{\Sob}_{s,\alpha,\bsgamma,q}$, is a Banach space consisting of non-periodic (in the sense of not necessarily periodic) smooth functions with the norm
\begin{eqnarray*}
 && \|f\|^{\Sob}_{s,\alpha,\bsgamma,q} := \\
 &&  \small
\sup_{u\subseteq \{1,\ldots,s\}}\gamma_u^{-1}\left(\sum_{v\subseteq u}\sum_{\bstau_{u\setminus v}\in \{1,\ldots,\alpha\}^{|u\setminus v|}}\int_{[0,1)^{|v|}}\left| \int_{[0,1)^{s-|v|}} f^{(\bstau_{u\setminus v},\bsalpha_v,\bszero)}(\bsx) \rd \bsx_{-v}\right|^{q}\rd \bsx_v\right)^{1/q}, 
\end{eqnarray*}
where $(\bstau_{u\setminus v},\bsalpha_v,\bszero)$ denotes the vector $\bsh\in \NN_0^s$ such that $h_j=\tau_j$ if $j\in u\setminus v$, $h_j=\alpha$ if $j\in v$, and $h_j=0$ otherwise, and $f^{(\bstau_{u\setminus v},\bsalpha_v,\bszero)}(\bsx)$ denotes the mixed derivative of order $(\bstau_{u\setminus v},\bsalpha_v,\bszero)$ of $f$. Moreover, we write $\bsx_v=(x_j)_{j\in v}$ and $\bsx_{-v}=(x_j)_{j\in \{1,\ldots,s\}\setminus v}$.
\end{definition}

This Sobolev space was introduced by \cite{DKLNS14} in the context of partial differential equations with random coefficients.
(The original function space in \cite{DKLNS14} contains the additional parameter $r\in [1,\infty]$, and the definition of the norm has been corrected in \url{https://arxiv.org/abs/1309.4624}. In this paper  we choose $r=\infty$, which makes the norm smallest over $r\in [1,\infty]$.) 
The parameter $\alpha$ determines the differentiability of the non-periodic functions. As for the rank-1 lattice rules for the weighted Korobov spaces, it is desirable to have good modulus $p$ and generating vector $\bsq$ such that the worst-case error of the corresponding high-order polynomial lattice rule for $\cF^{\Sob}_{s,\alpha,\bsgamma}$ is small. Originally in \cite{DKLNS14}, interlaced polynomial lattice rules \cite{G15,GD15} were used instead of high-order polynomial lattice rules, and it was shown that the worst-case error bound of order $N^{-\alpha+\epsilon}$ with arbitrarily small $\epsilon>0$ can be achieved by the CBC algorithm applied to interlaced polynomial lattice rules. The major advantage of interlaced polynomial lattice rules over high-order polynomial lattice rules lies in the construction cost for the CBC algorithm: for the product weights, constructing an interlaced rule requires $\cO(\alpha sN\log N)$ operations with $\cO(N)$ memory 
\cite{G15}, whereas constructing a high-order rule require $\cO(sN^{\alpha}\log N)$ operations with $\cO(N^{\alpha})$ memory \cite{BDLNP12}.

However, interlaced polynomial lattice rules are not necessarily a better choice than high-order polynomial lattice rules. To construct an interlaced polynomial lattice rule, which relies on the digit interlacing method due to Dick \cite{Dick08}, we must select an integer interlacing factor $d$, the construction cost increases linearly with $d$,
and the resulting rule cannot exploit the smoothness of functions beyond $d$.  This means that if $\alpha>d$, the worst-case error bound is only of order $N^{-d+\epsilon}$. High-order polynomial lattice rules do not explicitly require such a factor. We only need to specify the maximum precision $n$ of the points. This $n$ can be set as large as possible, for instance, with $b=2$, we can take $n=53$ for the double-precision floating-point format. 
This way, high-order polynomial lattice rules can be possibly made adaptive to the smoothness $\alpha$ of functions, addressing the drawback of interlaced polynomial lattice rules.
Since we do not apply CBC in this paper, we prefer high-order rules over the interlaced ones. 

In what follows, we assume that the polynomial $p$ is irreducible and we write
\[ 
  G_n := \{q\in \FF_b[x]\mid q\neq 0\quad \text{and}\quad \deg(q)<n\}.
\]
Analogously to the rank-1 lattice case in Section~\ref{sec:lattice}, 
we consider the following median high-order polynomial lattice rule for weighted Sobolev spaces. 
For an odd integer $r$, we draw $\bsq_1,\ldots,\bsq_r$ randomly and independently from the set $G_n^s$, 
and we approximate $I_s(f)$ by 
\[ 
   M_{m,n,s,p,r}(f)
	:= \median\left(Q_{P_{m,n,s,p,\bsq_1}}(f),\ldots,Q_{P_{m,n,s,p,\bsq_r}}(f)\right).
\]
The worst-case error is the random variable
\[
  e^{\wor}(M_{m,n,s,p,r}; \cF^{\Sob}_{s,\alpha,\bsgamma,q}) 
	:=  \sup_{\substack{f\in \cF^{\Sob}_{s,\alpha,\bsgamma,q}\\ 
	 \|f\|^{\Sob}_{s,\alpha,\bsgamma,q}\leq 1}} |M_{m,n,s,p,r}(f) -I_s(f)|. 
\]

\subsection{Main results for polynomial lattice point sets in Sobolev spaces}

We first need a few definitions and lemmas.

\begin{definition}[dual polynomial lattice]\label{def:dual_net}
Let $m,n\in \NN$ with $m\leq n$, $p\in \FF_b[x]$ with $\deg(p)=n$ and $\bsq \in G_n^s$. For $k\in \NN_0$ with the $b$-adic finite expansion $k=\kappa_0+\kappa_1b+\cdots$, we define
\[ \tr_n(k)=\sum_{i=0}^{n-1}\kappa_ix^i \in G_n\cup \{0\}.\]
This operator is applied component-wise to a vector. Then the set
\[ P^{\perp}_{m,n,s,p,\bsq} = \left\{ \bsk \in \NN_0^s\mid \tr_n(\bsk)\cdot \bsq\equiv a \pmod p \quad \text{with $\deg(a)<n-m$}\right\}, \]
is called the dual net of the high-order polynomial lattice point set $P_{m,n,s,p,\bsq}$.
\end{definition}

\begin{definition}[Walsh functions]\label{def:walsh}
Let us write $\omega_b:=\exp(2\pi \ri/b)$. For $k\in \NN_0$, we denote the $b$-adic expansion of $k$ by $k=\kappa_0+\kappa_1b+\cdots$. The $k$-th Walsh function $\wal_k\colon [0,1)\to \CC$ is defined by
\[ \wal_k(x) := \omega_b^{\kappa_0\xi_1+\kappa_1\xi_2+\cdots}, \]
where the $b$-adic expansion of $x\in [0,1)$ is denoted by $x=\xi_1/b+\xi_2/b^2+\cdots$, which is understood to be unique in the sense that infinitely many of the $\xi_i$ are different from $b-1$.

For $s\geq 2$ and $\bsk=(k_1,\ldots,k_s)\in \NN_0^s$, the $s$-dimensional $\bsk$-th Walsh function $\wal_{\bsk} \colon [0,1)^s\to \CC$ is defined by
\[ \wal_{\bsk}(\bsx) := \prod_{j=1}^{s}\wal_{k_j}(x_j). \]
\end{definition}
\noindent It is well-known that the system of
Walsh functions is a complete orthogonal system in $L_2([0,1)^s)$, see \cite[Appendix~A]{DickPilli10}. The following character property of the high-order polynomial lattice point set is analogous to what is stated in Lemma~\ref{lem:character}.
\begin{lemma}[character property]\label{lem:character_2}
For $m,n\in \NN$ with $m\leq n$, $p\in \FF_b[x]$ with $\deg(p)=n$ and $\bsq \in G_n^s$, we have
\[ \frac{1}{b^m}\sum_{\bsx\in P_{m,n,s,p,\bsq}}\wal_{\bsk}(\bsx)=\begin{cases} 1 & \text{if $\bsk\in P^{\perp}_{m,n,s,p,\bsq}$,} \\ 0 & \text{otherwise.} \end{cases}\]
\end{lemma}

For any $f\in \cF^{\Sob}_{s,\alpha,\bsgamma,q}$, we have the following absolutely convergent Walsh series
\[ f(\bsx) = \sum_{\bsk\in \NN_0^s}\tilde{f}(\bsk)\wal_{\bsk}(\bsx),\]
where $\tilde{f}(\bsk)$ denotes the $\bsk$-th Walsh coefficient of $f$:
\[ \tilde{f}(\bsk) := \int_{[0,1)^s}f(\bsx)\overline{\wal_{\bsk}(\bsx)}\rd \bsx. \]
Note that $\tilde{f}(\bszero)$ coincides with the integral $I_s(f)$. The following result on the decay of Walsh coefficients for $f\in \cF^{\Sob}_{s,\alpha,\bsgamma,q}$ was shown in \cite[Theorem~3.5]{DKLNS14}.

\begin{lemma}[decay of Walsh coefficients]\label{lem:walsh_decay}
Let $\alpha\in \NN$, $\alpha\geq 2$, $1\leq q\leq \infty$ and $\bsgamma=\{\gamma_u\}_{u\subseteq \{1,\ldots,s\}}$ be a set of positive weights with $\gamma_{\emptyset}=1$. For any $f\in \cF^{\Sob}_{s,\alpha,\bsgamma,q}$, a non-empty subset $u\subseteq \{1,\ldots,s\}$ and $\bsk_u\in \NN^{|u|}$, it holds that
\[ |\tilde{f}(\bsk_u,\bszero)|\leq \|f\|^{\Sob}_{s,\alpha,\bsgamma,q}\gamma_uC_{\alpha}^{|u|}b^{-\mu_{\alpha}(\bsk_u)},\]
where
\[\small 
  C_{\alpha}:= \left(1+\frac{1}{b}+\frac{1}{b(b+1)} \right)^{\alpha-2}\left(3+\frac{2}{b}+\frac{2b+1}{b-1} \right)\max\left( \frac{2}{(2\sin\frac{\pi}{b})^{\alpha}},\max_{1\leq \tau< \alpha}\frac{1}{(2\sin\frac{\pi}{b})^{\tau}}\right)
\]
and $\mu_{\alpha}(\bsk_u):=\sum_{j\in u}\mu_{\alpha}(k_j)$, with
\[ \mu_{\alpha}(k)=\sum_{i=1}^{\min(\alpha,c)}a_i\]
for $k\in \NN$ whose $b$-adic expansion is given by $k=\kappa_1b^{a_1-1}+\kappa_2b^{a_2-1}+\cdots+\kappa_cb^{a_c-1}$ such that $c\geq 1$, $\kappa_1,\ldots,\kappa_c\in \{1,\ldots,b-1\}$ and $a_1>\cdots>a_c$.
\end{lemma}

As the second main result of this paper, we show a probabilistic upper bound on the worst-case error of our median high-order polynomial lattice rule for weighted Sobolev spaces.

\begin{theorem}\label{thm:main2}
Let $m,n\in \NN$ with $m\leq n$, $p\in \FF_b[x]$ be irreducible with $\deg(p)=n$, $r$ be odd and $\bsq_1,\ldots,\bsq_r$ be chosen independently and randomly from the set $G_n^s$. Then, for any integer $\alpha\geq 2$ and $\bsgamma$, the worst-case error is bounded above by
\[\small
 e^{\wor}(M_{m,n,s,p,r}; \cF^{\Sob}_{s,\alpha,\bsgamma,q}) 
 \leq \inf_{1/\alpha<\lambda<1}\left(\frac{2}{\eta(b^{\min(m,\lambda n)}-1)}\sum_{\emptyset \neq u\subseteq \{1,\ldots,s\}}\gamma_u^{\lambda}C_{\alpha}^{\lambda|u|}A_{\alpha,\lambda}^{|u|}\right)^{1/\lambda}
\]
with a probability of at least
\[ 1-\binom{r}{(r+1)/2}\eta^{(r+1)/2},\]
for any $0<\eta<1$, where we write
\begin{align}\label{eq:A_alpha_lambda}
  A_{\alpha,\lambda}=\sum_{\tau=1}^{\alpha-1}\prod_{i=1}^{\tau}\frac{b-1}{b^{\lambda i}-1}+\frac{b^{\lambda \alpha}-1}{b^{\lambda \alpha}-b}\prod_{i=1}^{\alpha}\frac{b-1}{b^{\lambda i}-1}.
\end{align}
\end{theorem}

\begin{proof}
Throughout this proof, we write
\[ 
 \tilde{r}_{\alpha,\bsgamma}(\bsk_u,\bszero) = \gamma_u C_{\alpha}^{|u|} b^{-\mu_{\alpha}(\bsk_u)},
\]
for a non-empty subset $u\subseteq \{1,\ldots,s\}$ and $\bsk_u\in \NN^{|u|}$.
Since any $f\in \cF^{\Sob}_{s,\alpha,\bsgamma,q}$ has an absolutely convergent Walsh series, by applying Lemma~\ref{lem:character_2}, Lemma~\ref{lem:median_jensen},  H\"{o}lder's inequality and Lemma~\ref{lem:walsh_decay} in this order, it holds that
\begin{align}\label{eq:proof_main2_1}
    & \ \ \ \ \ \ e^{\wor}(M_{m,n,s,p,r}; \cF^{\Sob}_{s,\alpha,\bsgamma,q}) \notag\\
		& = \sup_{\substack{f\in \cF^{\Sob}_{s,\alpha,\bsgamma,q} \notag \\ \|f\|^{\Sob}_{s,\alpha,\bsgamma,q}\leq 1}}\left| \median_{1\leq \ell\leq r}\frac{1}{b^m}\sum_{\bsx\in P_{m,n,s,p,\bsq_{\ell}}}f(\bsx)-I(f)\right| \notag \\
    & = \sup_{\substack{f\in \cF^{\Sob}_{s,\alpha,\bsgamma,q} \notag \\ \|f\|^{\Sob}_{s,\alpha,\bsgamma,q}\leq 1}}\left| \median_{1\leq \ell\leq r}\frac{1}{b^m}\sum_{\bsx\in P_{m,n,s,p,\bsq_{\ell}}}\sum_{\bsk\in \NN_0^s}\tilde{f}(\bsk)\wal_{\bsk}(\bsx)-\tilde{f}(\bszero)\right| \notag \\
    & = \sup_{\substack{f\in \cF^{\Sob}_{s,\alpha,\bsgamma,q} \notag \\ \|f\|^{\Sob}_{s,\alpha,\bsgamma,q}\leq 1}}\left| \median_{1\leq \ell\leq r}\sum_{\bsk\in P^{\perp}_{m,n,s,p,\bsq_{\ell}}\setminus \{\bszero\}}\tilde{f}(\bsk)\right| \notag \\
    & \leq \sup_{\substack{f\in \cF^{\Sob}_{s,\alpha,\bsgamma,q} \notag \\ \|f\|^{\Sob}_{s,\alpha,\bsgamma,q}\leq 1}} \median_{1\leq \ell\leq r}\sum_{\bsk\in P^{\perp}_{m,n,s,p,\bsq_{\ell}}\setminus \{\bszero\}}\frac{|\tilde{f}(\bsk)|}{\tilde{r}_{\alpha,\bsgamma}(\bsk)}\tilde{r}_{\alpha,\bsgamma}(\bsk) \notag \\
    & \leq \sup_{\substack{f\in \cF^{\Sob}_{s,\alpha,\bsgamma,q} \notag \\ 
\|f\|^{\Sob}_{s,\alpha,\bsgamma,q}\leq 1}} \median_{1\leq \ell\leq r}\left(\sup_{\bsk\in P^{\perp}_{m,n,s,p,\bsq_{\ell}}\setminus \{\bszero\}}\frac{|\tilde{f}(\bsk)|}{\tilde{r}_{\alpha,\bsgamma}(\bsk)}\right)\left(\sum_{\bsk\in P^{\perp}_{m,n,s,p,\bsq_{\ell}}\setminus \{\bszero\}}\tilde{r}_{\alpha,\bsgamma}(\bsk)\right) \notag \\
    & \leq \sup_{\substack{f\in \cF^{\Sob}_{s,\alpha,\bsgamma,q} \notag \\ \|f\|^{\Sob}_{s,\alpha,\bsgamma,q}\leq 1}} \left(\sup_{\bsk\in \NN_0^s\setminus \{\bszero\}}\frac{|\tilde{f}(\bsk)|}{\tilde{r}_{\alpha,\bsgamma}(\bsk)}\right) 
		\median_{1\leq \ell\leq r}  \tilde{\Scal}_{\alpha,\bsgamma,p}(\bsq_\ell) \notag \\
    & \leq \median_{1\leq \ell\leq r}  \tilde{\Scal}_{\alpha,\bsgamma,p}(\bsq_\ell),  
\end{align}
where
\[ 
  \tilde{\Scal}_{\alpha,\bsgamma,p}(\bsq)
	= \sum_{\bsk\in P^{\perp}_{m,n,s,p,\bsq}\setminus \{\bszero\}}\tilde{r}_{\alpha,\bsgamma}(\bsk). 
\]
For $1/(\alpha)<\lambda\leq 1$, by using the subadditivity \eqref{eq:jensen}, we have
\begin{align*}
    \frac{1}{|G_n|^s}\sum_{\bsq\in G_n^s}(\tilde{\Scal}_{\alpha,\bsgamma,p}(\bsq))^{\lambda} & = \frac{1}{(b^n-1)^s}\sum_{\bsq\in G_n^s}\left(\sum_{\bsk\in P^{\perp}_{m,n,s,p,\bsq_{\ell}}\setminus \{\bszero\}}\tilde{r}_{\alpha,\bsgamma}(\bsk)\right)^{\lambda} \\
    & \leq \frac{1}{(b^n-1)^s}\sum_{\bsq\in G_n^s}\sum_{\bsk\in P^{\perp}_{m,n,s,p,\bsq_{\ell}}\setminus \{\bszero\}}(\tilde{r}_{\alpha,\bsgamma}(\bsk))^{\lambda}\\
    & = \sum_{\bsk\in \NN_0^s\setminus \{\bszero\}}(\tilde{r}_{\alpha,\bsgamma}(\bsk))^{\lambda}\frac{1}{(b^n-1)^s}\sum_{\substack{\bsq\in G_n^s\\ \tr_n(\bsk)\cdot \bsq \equiv a\pmod p\\ \deg(a)<n-m}}1.
\end{align*}
If $p\mid \tr_n(\bsk)$, the condition $\tr_n(\bsk)\cdot \bsq \equiv a\pmod p$ trivially holds with $a=0$ for all $\bsq\in G_n^s$. Otherwise if $p\nmid \tr_n(\bsk)$, i.e., if there exists a non-empty subset $u\subseteq \{1,\ldots,s\}$ such that $p\nmid \tr_n(k_j)$ for all $j\in u$ and $p\mid \tr_n(k_j)$ for $j\not\in u$, the condition $\tr_n(\bsk)\cdot \bsq \equiv a\pmod p$ is equivalent to $\tr_n(\bsk_u)\cdot \bsq_u \equiv a\pmod p$, which itself is equivalent to
\[ \tr_n(k_j) q_j \equiv a-\tr_n(\bsk_{u\setminus \{j\}})\cdot \bsq_{u\setminus \{j\}}\pmod p, \]
for any $j\in u$. As we have $p\nmid \tr_n(k_j)$ and we assume that $p$ is irreducible, there exists at most one $q_j\in G_n$ which satisfies the above equality for each $a\in \FF_b[x]$ with $\deg(a)<n-m$ and $\bsq_{u\setminus \{j\}}\in G_n^{|u|-1}$. Therefore, the number of $\bsq\in G_n^s$ which satisfy $p\nmid \tr_n(\bsk)$ and $\tr_n(\bsk)\cdot \bsq \equiv a\pmod p$ with $\deg(a)<n-m$ is bounded above by the product of the number of possible choices for $a\in \FF_b[x]$, which is $b^{n-m}$, and the number of possible choices for $\bsq_{\{1,\ldots,s\}\setminus \{j\}}$, which is $(b^n-1)^{s-1}$. Thus it follows that
\begin{align*}
    & \frac{1}{|G_n|^s}\sum_{\bsq\in G_n^s}(\tilde{\Scal}_{\alpha,\bsgamma,p}(\bsq))^{\lambda} \\ 
    & \leq \sum_{\substack{\bsk\in \NN_0^s\setminus \{\bszero\}\\ p\mid \tr_n(\bsk)}}(\tilde{r}_{\alpha,\bsgamma}(\bsk))^{\lambda}+\sum_{\substack{\bsk\in \NN_0^s\setminus \{\bszero\}\\ p\nmid \tr_n(\bsk)}}(\tilde{r}_{\alpha,\bsgamma}(\bsk))^{\lambda}\frac{b^{n-m}(b^n-1)^{s-1}}{(b^n-1)^s} \\
    & \leq \sum_{\bsk\in \NN_0^s\setminus \{\bszero\}}(\tilde{r}_{\alpha,\bsgamma}(b^n\bsk))^{\lambda}+\frac{1}{b^m-1}\sum_{\bsk\in \NN_0^s\setminus \{\bszero\}}(\tilde{r}_{\alpha,\bsgamma}(\bsk))^{\lambda}\\
    & = \sum_{\emptyset \neq u\subseteq \{1,\ldots,s\}}\sum_{\bsk_u\in \NN^{|u|}}(\tilde{r}_{\alpha,\bsgamma}(b^n\bsk_u,\bszero))^{\lambda}+\frac{1}{b^m-1}\sum_{\emptyset \neq u\subseteq \{1,\ldots,s\}}\sum_{\bsk_u\in \NN^{|u|}}(\tilde{r}_{\alpha,\bsgamma}(\bsk_u,\bszero))^{\lambda}\\
    & = \sum_{\emptyset \neq u\subseteq \{1,\ldots,s\}}\gamma_u^{\lambda}C_{\alpha}^{\lambda|u|}\left(\sum_{k\in \NN}b^{-\lambda\mu_{\alpha}(b^nk)}\right)^{|u|} \\
		& \ \ \ \ \ + \frac{1}{b^m-1}\sum_{\emptyset \neq u\subseteq \{1,\ldots,s\}}\gamma_u^{\lambda}C_{\alpha}^{\lambda|u|}\left(\sum_{k\in \NN}b^{-\lambda\mu_{\alpha}(k)}\right)^{|u|}\\
    & \leq \sum_{\emptyset \neq u\subseteq \{1,\ldots,s\}}\gamma_u^{\lambda}\frac{C_{\alpha}^{\lambda|u|}A_{\alpha,\lambda}^{|u|}}{b^{\lambda n|u|}}+\frac{1}{b^m-1}\sum_{\emptyset \neq u\subseteq \{1,\ldots,s\}}\gamma_u^{\lambda}C_{\alpha}^{\lambda|u|}A_{\alpha,\lambda}^{|u|}\\
    & \leq \frac{2}{b^{\min(m,\lambda n)}-1}\sum_{\emptyset \neq u\subseteq \{1,\ldots,s\}}\gamma_u^{\lambda}C_{\alpha}^{\lambda|u|}A_{\alpha,\lambda}^{|u|},
\end{align*}
where we have used the results of \cite[Lemma~7]{G16} on the sums of $b^{-\lambda\mu_{\alpha}(k)}$ and $b^{-\lambda\mu_{\alpha}(b^nk)}$ in the third inequality, which involve $A_{\alpha,\lambda}$ given in \eqref{eq:A_alpha_lambda}. This gives a bound on the average of $(\tilde{\Scal}_{\alpha,\bsgamma,p}(\bsq))^{\lambda}$ which holds for any $1/\alpha<\lambda\leq 1$. 

Then, Markov's inequality ensures that, for any $0<\eta<1$, the event 
\[ 
  \tilde{\Scal}_{\alpha,\bsgamma,p}(\bsq) > \inf_{1/\alpha<\lambda<1}\left(\frac{2}{\eta(b^{\min(m,\lambda n)}-1)}\sum_{\emptyset \neq u\subseteq \{1,\ldots,s\}}\gamma_u^{\lambda}C_{\alpha}^{\lambda|u|}A_{\alpha,\lambda}^{|u|}\right)^{1/\lambda}
 =: \tilde B(\alpha,\bsgamma)
\]
happens with a probability of at most $\eta$ under a random choice of $\bsq\in G_n^s$.  
For the median estimator $M_{m,n,s,p,r}$ to be larger than this bound $\tilde B(\alpha,\bsgamma)$,
we must have $\tilde{\Scal}_{\alpha,\bsgamma}(\bsq_\ell) > \tilde B(\alpha,\bsgamma)$ for at least 
$(r+1)/2$ vectors among $\bsq_1,\ldots,\bsq_r$.
The probability that this happens is bounded above by
\[ 
  \binom{r}{(r+1)/2}\eta^{(r+1)/2}.
\]
Combining this with the bound shown in \eqref{eq:proof_main2_1} completes the proof.
\end{proof}

As pointed out in \cite[Section~3.1]{DKLNS14}, $\alpha\geq 2$ is required to ensure the convergence of the infinite sum 
\[ \sum_{\bsk\in P^{\perp}_{m,n,s,p,\bsq}\setminus \{\bszero\}}b^{-\mu_{\alpha}(\bsk)},\]
for any irreducible $p$ and $\bsq\in G_n^s$. Thus, the case $\alpha=1$ is not covered by our result.

Using Remark~\ref{rem:probability}, we obtain the following corollary:

\begin{corollary}
Let $\alpha\ge 2$, $\bsgamma$ be a set of weights, and $n\geq \alpha m = \alpha \log_b N$.
Then for any odd $r\ge 3$, $\epsilon > 0$ and $0 < \rho < 1$, 
there is a constant $c_1 = c_1(\alpha,\bsgamma,\epsilon) > 0$ 
(which depends on $\alpha$, the $\gamma_u$'s and $\epsilon$) such that
\[
  \PP\left[e^{\wor}(M_{m,n,s,p,r}; \cF^{\Sob}_{s,\alpha,\bsgamma})
	         \le \frac{c_1(\alpha,\bsgamma,\epsilon)}{(\rho N)^{\alpha-\epsilon}}\right] 
  \ge 1 - \rho^{(r+1)/2}/4.
\]
\end{corollary}

\begin{proof}
Take $1/\lambda = \alpha-\epsilon$.  Under the assumption on $n$, we have $n \ge \alpha m > m/\lambda$ 
and then $(b^{\min(m, \lambda n)})^{-1/\lambda} = b^{-\min(m/\lambda, n)} = N^{-1/\lambda} = N^{-\alpha + \epsilon}$. Then the result follows from Theorem~\ref{thm:main2} and 
the bound (\ref{eq:binom-bound}) in Remark~\ref{rem:probability}.
\end{proof}

Thus, provided that we take $n$ large enough, we get a convergence rate of almost $\cO(N^{-\alpha})$
(with high probability) for any $\alpha \ge 2$.  
In other words, our median high-order polynomial lattice rule exploits the smoothness of functions adaptively.
Note that Remark~\ref{rem:binomial} also applies here.

\section{Numerical experiments}
\label{sec:experiments}

We conclude this paper with numerical experiments both for rank-1 lattice rules and high-order polynomial lattice rules.
The goal is to illustrate how the worst-case error for the median rule truly behaves on some concrete examples.
In particular, we want to illustrate the fact that most of the possible generating vectors are a good choice, 
while a small minority are bad.

\subsection{Lattice rules for periodic functions}

\begin{figure}
    \centering
    \includegraphics[width=0.45\linewidth]{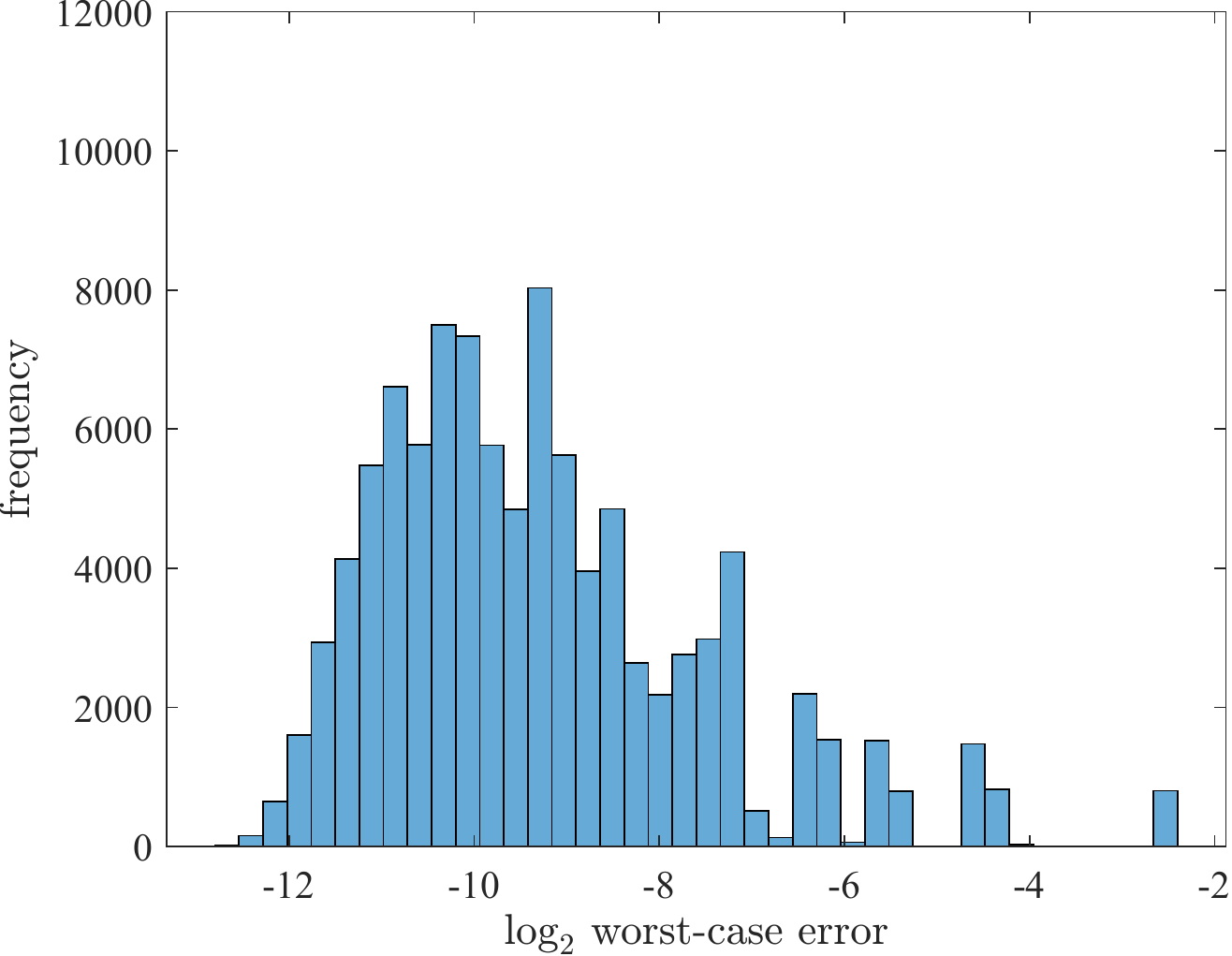}
    \includegraphics[width=0.45\linewidth]{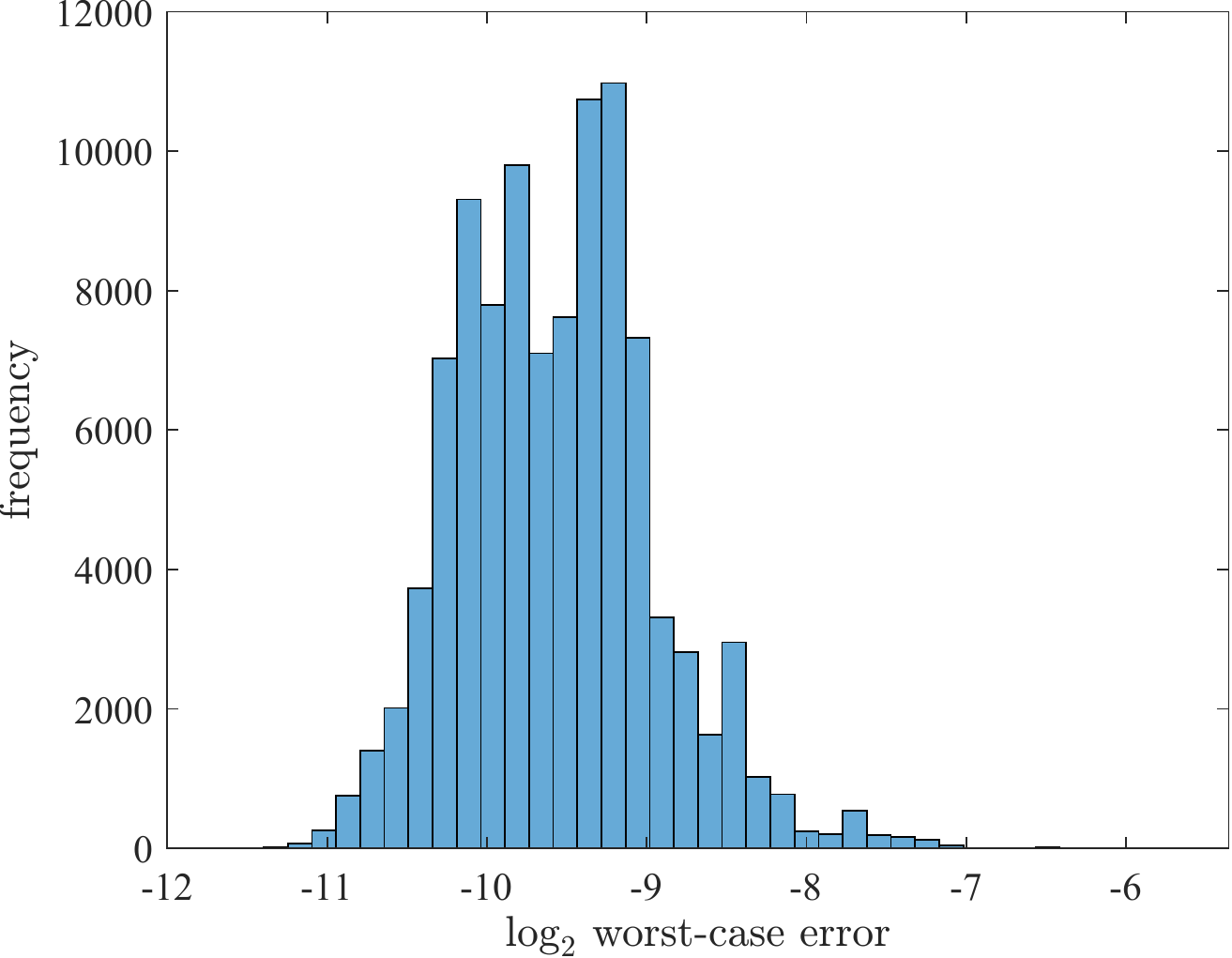}\\[10pt]
    \includegraphics[width=0.45\linewidth]{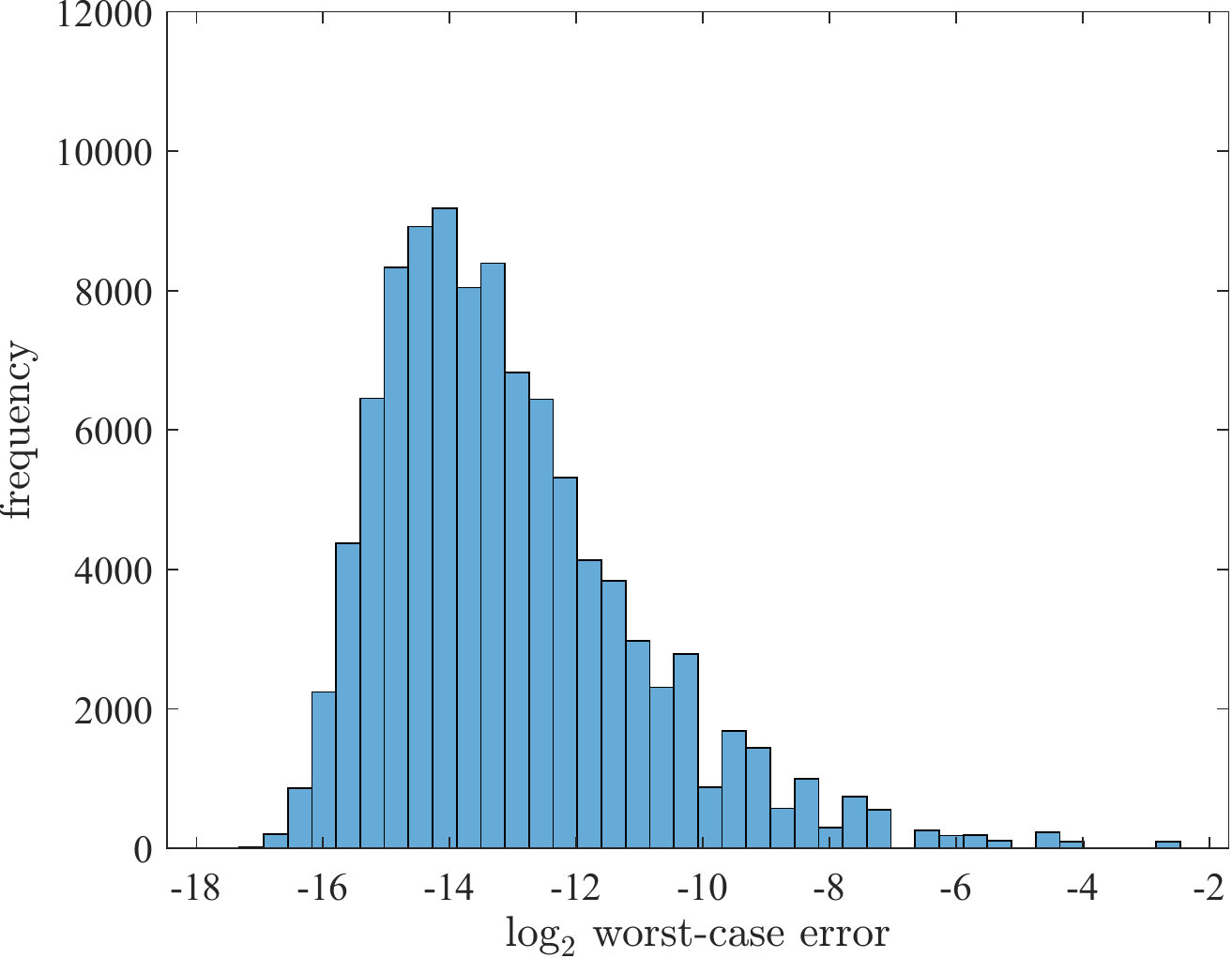}
    \includegraphics[width=0.45\linewidth]{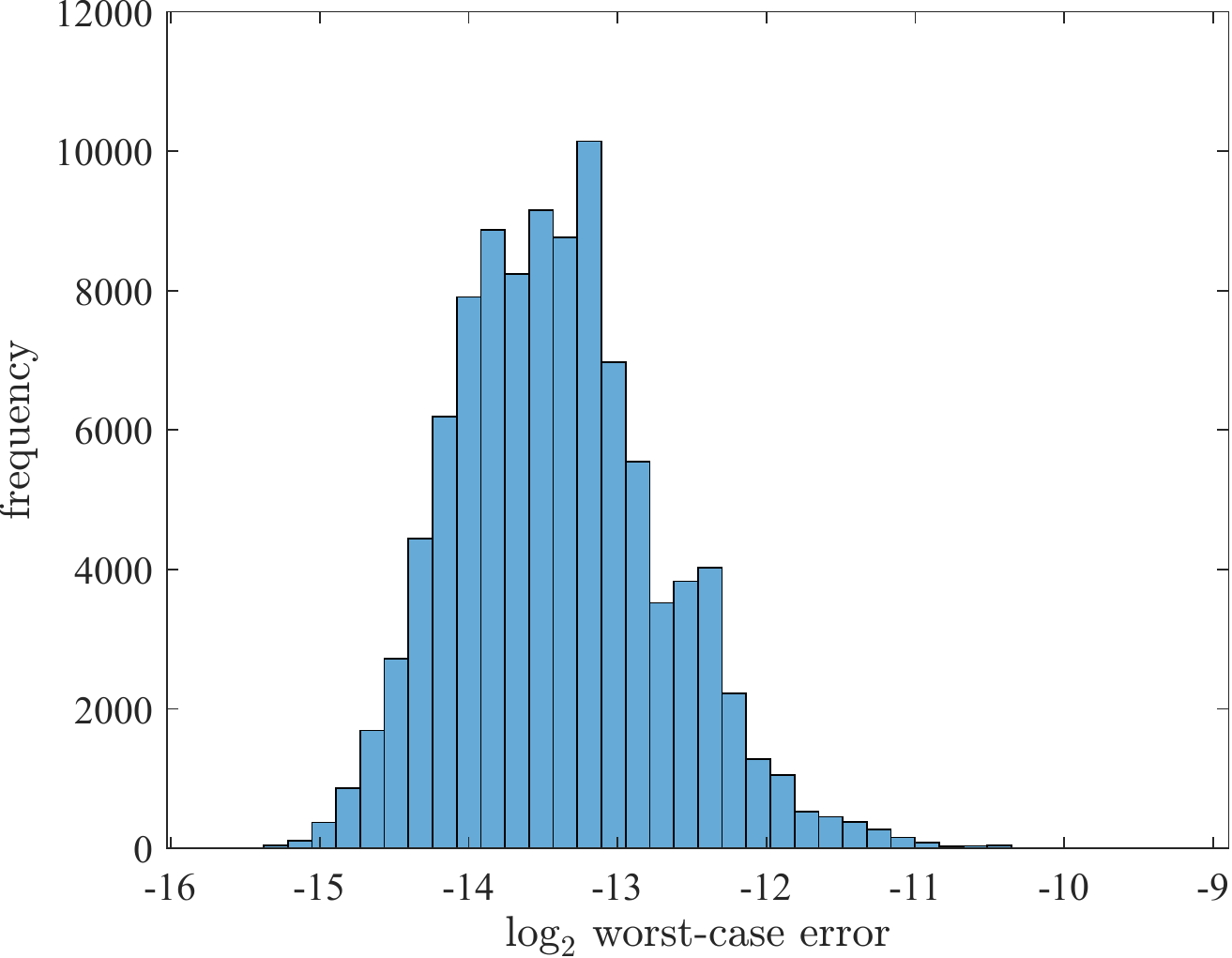}\\
    \caption{Histograms of the $\log_2$ of the worst-case error $\Scal_{\alpha,\bsgamma}(\bsz)$ with $\alpha=2$ and $\gamma_j=1/j^3$ for rank-1 lattice rules with randomly chosen generating vectors with $N=251$ (upper panels) and $N=2039$ (lower panels). The left panels are for a single choice ($r=1$), while for the right panels we take the median of the worst-case error for rank-1 lattice rules with $r=11$ randomly chosen generating vectors.}
    \label{fig:dist_wce}
\end{figure}

\begin{example}\rm 
For our first example, we consider a weighted Korobov space with integer smoothness parameter $\alpha \ge 1$ and 
product weights $\gamma_u = \prod_{j\in u} \gamma_j$.  The worst-case error of the rank-1 lattice rule with generating vector $\bsz$ for that space has the explicit form
\[ 
  \Scal_{\alpha,\bsgamma}(\bsz) 
 = \left(-1+\frac{1}{N}\sum_{\bsx\in P_{N,s,\bsz}}
	 \prod_{j=1}^{s}\left[1+\gamma^2_j\frac{(-1)^{\alpha+1}(2\pi)^{2\alpha}}{(2\alpha)!}B_{2\alpha}(x_j) \right]\right)^{1/2}, 
\]
where $B_{2\alpha}$ denotes the Bernoulli polynomial of degree $2\alpha$; see \cite{vLEC12a} and \cite[Section~5]{DKS13}. 
In this artificial simple case, we know the exact optimal weights that must be taken in a CBC search for $\bsz$,
so we can compare the median estimator with the best possible case of a CBC search.

We take two primes $N=251$ and $N=2039$, both for $s=50$ dimensions, with $\alpha=2$ and $\gamma_j=1/j^3$. 
For each of those $N$, we drew $10^5$ generating vectors $\bsz$ randomly and uniformly from $\{1,\ldots,N-1\}^s$, 
and computed $\Scal_{\alpha,\bsgamma}(\bsz)$ for each.  
The left panels of Figure~\ref{fig:dist_wce} show a histogram of the $10^5$ realizations of 
$\log_2 \Scal_{\alpha,\bsgamma}(\bsz)$ for each of these two cases. 
Each histogram provides a good estimate of the true distribution of $\log_2 \Scal_{\alpha,\bsgamma}(\bsz)$,
which is a discrete distribution because $\bsz$ is drawn from a finite set.
Interestingly, the distributions are very asymmetric and are far from smooth on the right side: 
some rectangles are very high while others are zero in the same area. 
The largest observed values are $-2.4353$ for $N=251$ and $-2.4967$ for $N=2039$.
We can estimate from this data the $q$-quantiles $y_q$ of the distribution of 
$\Scal_{\alpha,\bsgamma}(\bsz)$, similar to those of the distribution of $Q_{P_{N,s,\bsz}}(f)$ 
in Remark~\ref{rem:binomial}.
For $q=0.75$, the corresponding empirical $q$-quantiles are $-8.3907$ for $N=251$ and $-12.0306$ for $N=2039$, 
while, for $q=0.9$, they are $-7.0975$ for $N=251$ and $-10.3101$ for $N=2039$.
These quantiles are much less than the worst observed values.
These empirical results agree with the fact that only a very small proportion of the vectors $\bsz$ are bad.
Suppose we draw $r$ random realizations of $\bsz$ and want the median of the $r$ corresponding 
values of $\Scal_{\alpha,\bsgamma}(\bsz)$ to be larger than 
$y = 10^{-3} \approx 2^{-10}$ with a probability smaller than $10^{-4}$.
For $N=2039$ this $y$ equals $y_q$ for $q \approx 0.9$, and Figure~\ref{fig:quantile_prob} shows that we can 
achieve approximately the target probability of $10^{-4}$ with $r=13$. 
For a larger $N$, the required $r$ is even smaller.
Note that for $N = 251$, $y_p = 2^{-10}$ corresponds to some $q < 0.5$, for which
the target probability of $10^{-4}$ cannot be achieved even for a very large $r$, 
as shown in Figure~\ref{fig:quantile_prob}.

For the remainder of our experiments reported in this paper, we took $r=11$.
The right panels of Figure~\ref{fig:dist_wce} show histograms of $10^5$ independent realizations of
$\log_2 [\median (\Scal_{\alpha,\bsgamma}(\bsz_{1}),\? \dots, \?\Scal_{\alpha,\bsgamma}(\bsz_{r}))]$ 
for randomly chosen $\bsz_1,\ldots,\bsz_r$ with $r=11$, corresponding to the cases $N=251$ and $N=2039$.
We see that the distributions have much less variance and are more symmetric than for a single random $\bsz$,
confirming the fact that taking the median successfully filters (adaptively) the bad vector generators.
Recall that the standard deviation of the empirical median as a function of $r$ 
generally decreases as $\cO(r^{-1/2})$.  That is, increasing $r$ decreases the noise rather slowly.
For the following examples, we made additional experiments with $r=31$ to see if it would make the 
error plots less noisy, and we did not see much visible difference.
\end{example}

\begin{example}\rm 
For our second example, we perform a numerical integration of the smooth, periodic function
\[ 
  f^{\per}_{\beta,\bsomega}(\bsx)
	 =\prod_{j=1}^{s}\left[ 1+\omega_j \left( g_{\beta}(x_j)-1\right)\right],
\]
with parameters $\beta$ and $\omega_j$, where the univariate function $g_{\beta}: [0,1]\to \RR$ is defined by
\[ 
 g_{\beta}(x)= (2\beta+1)\binom{2\beta}{\beta}x^{\beta}(1-x)^{\beta}.
\]
Note that $I_s(f^{\per}_{\beta,\bsomega})=1$. 
The function $g_{\beta}$ has been used for periodization of non-periodic functions, and our test function $f^{\per}_{\beta,\bsomega}$ belongs to the Korobov space with $\alpha=\beta$ when $\beta$ is a positive integer, see \cite[Section~5.10]{DKS13}. In what follows, we take $s=50$ and consider the four cases that correspond to $\beta=2$ or $\beta=5$, and  $\omega_j=1/j^{\beta+1}$ or $\omega_j=1/(s-j+1)^{\beta+1}$. 
We compare our median lattice rule with $r=11$, a QMC rule using non-randomized Sobol' points provided by MATLAB, 
and the rank-1 lattice rule with generating vector constructed by the fast CBC algorithm with 
$\Scal_{\alpha,\bsgamma}(\bsz)$ as a criterion, with $\alpha=2$ and the product weights $\gamma_j = 1/j^{3}$. 
These weights are not optimal, but they are a good heuristic choice when $\omega_j=1/j^{\beta+1}$.
When $\omega_j=1/(s-j+1)^{\beta+1}$, on the other hand, the weights decrease in the opposite 
direction as they should: they are very large for the unimportant coordinates and small for the important ones.
We do this to show how badly the CBC construction method can work when we have the wrong weights,
whereas the median estimator does not need any knowledge about the weights to perform well.
We choose $N$ to be a power of 2 for Sobol' points and to be a prime close to a power of 2 for lattice point sets.

The results for the four cases are shown in the corresponding panels of Figure~\ref{fig:lattice}. Both our median lattice rule and the rank-1 lattice rule constructed by the CBC algorithm can exploit the periodicity of the integrand and achieve a higher-order rate of convergence than $\cO(1/N)$. With a good choice of the weights in the CBC algorithm, the resulting rank-1 lattice rule performs better than our median lattice rule, as shown in the left panels. However, as the right panels clearly depict, if the relative importance of each of individual variables is not correctly specified, the performance of the rank-1 lattice rule with the CBC algorithm can deteriorate and even become inferior to the QMC rule using the Sobol' points when $N$ is not large. In contrast, our median lattice rule performs quite stably regardless of smoothness and weights.

\begin{figure}
    \centering
    \includegraphics[width=0.45\linewidth]{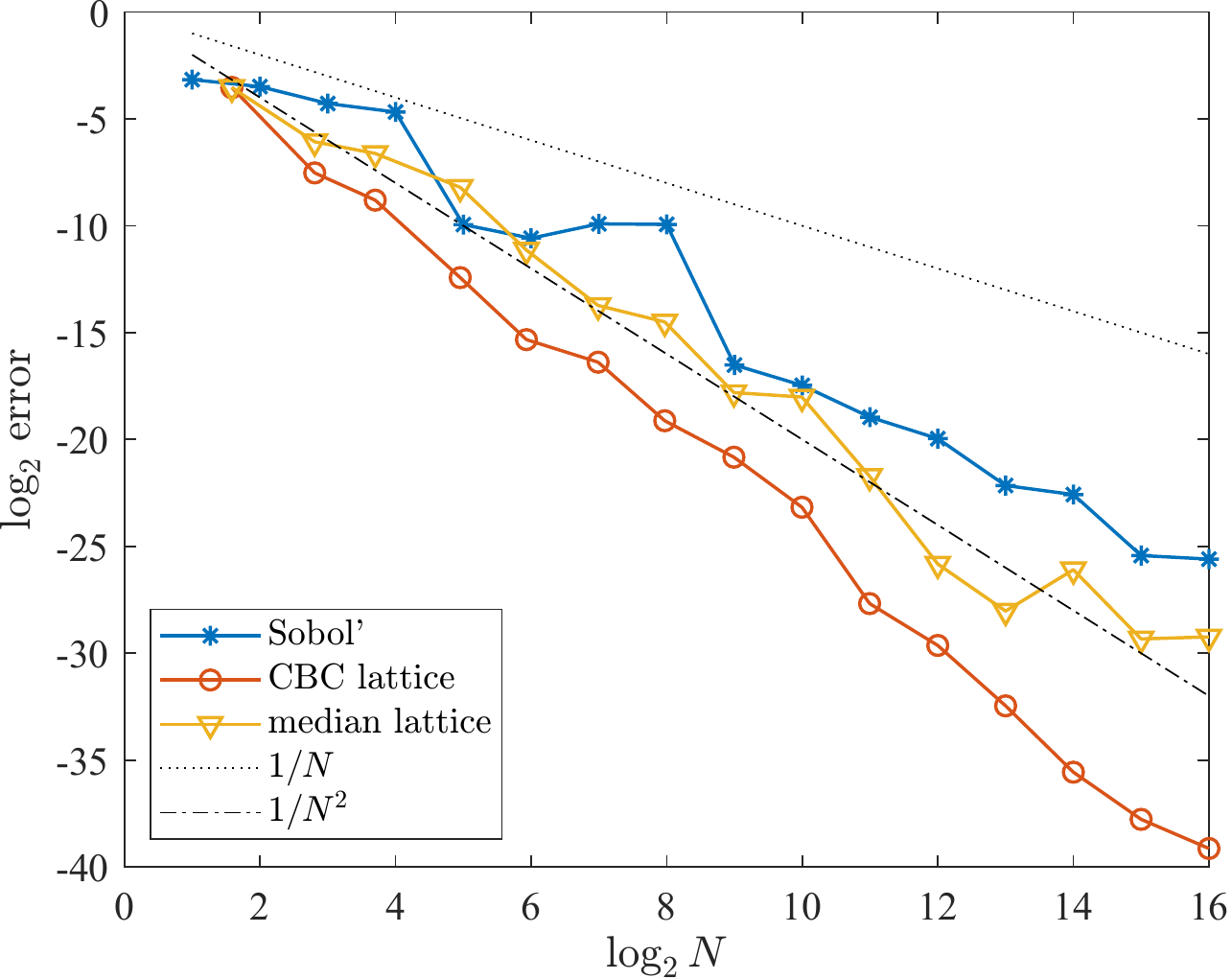}
    \includegraphics[width=0.45\linewidth]{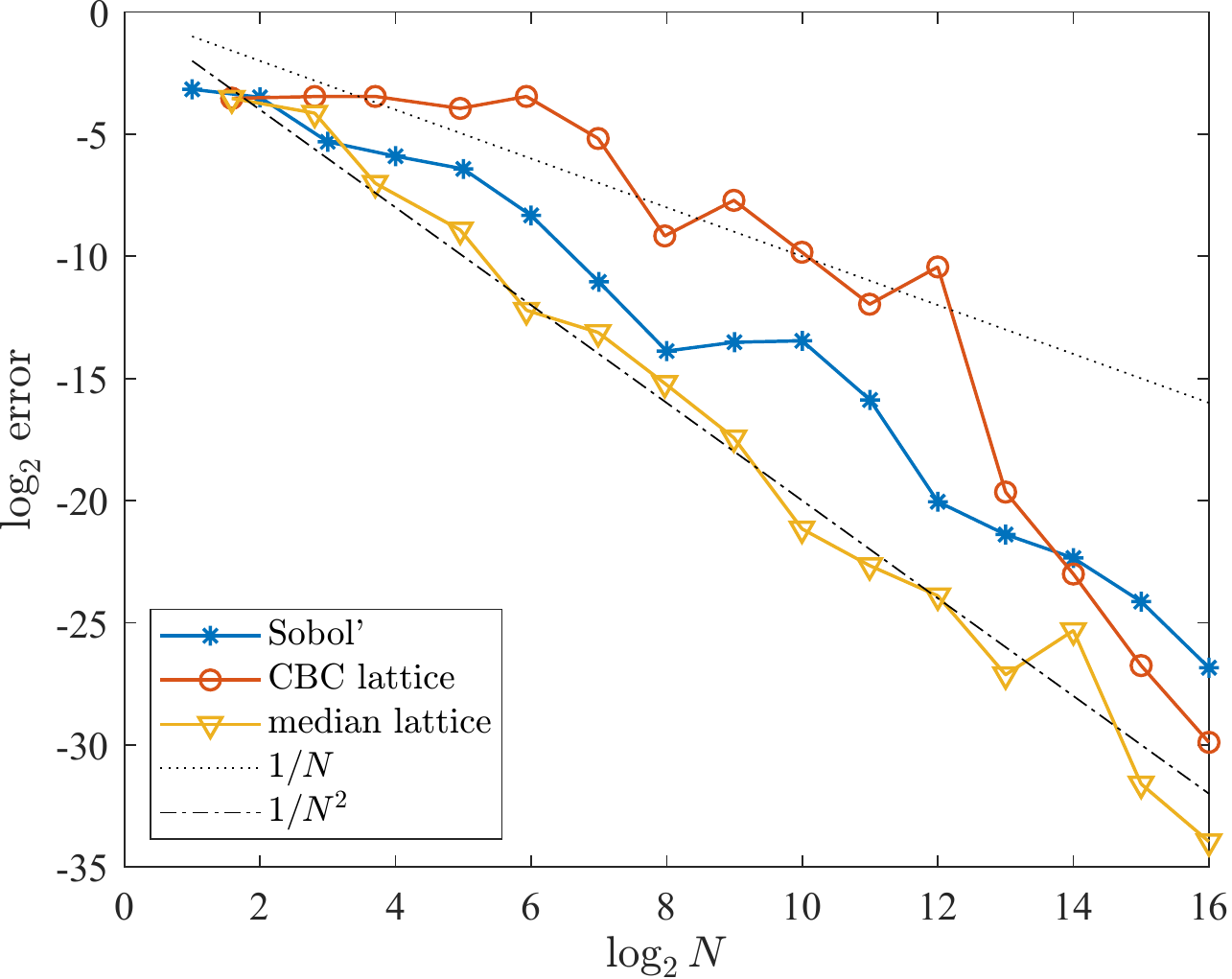}\\[10pt]
    \includegraphics[width=0.45\linewidth]{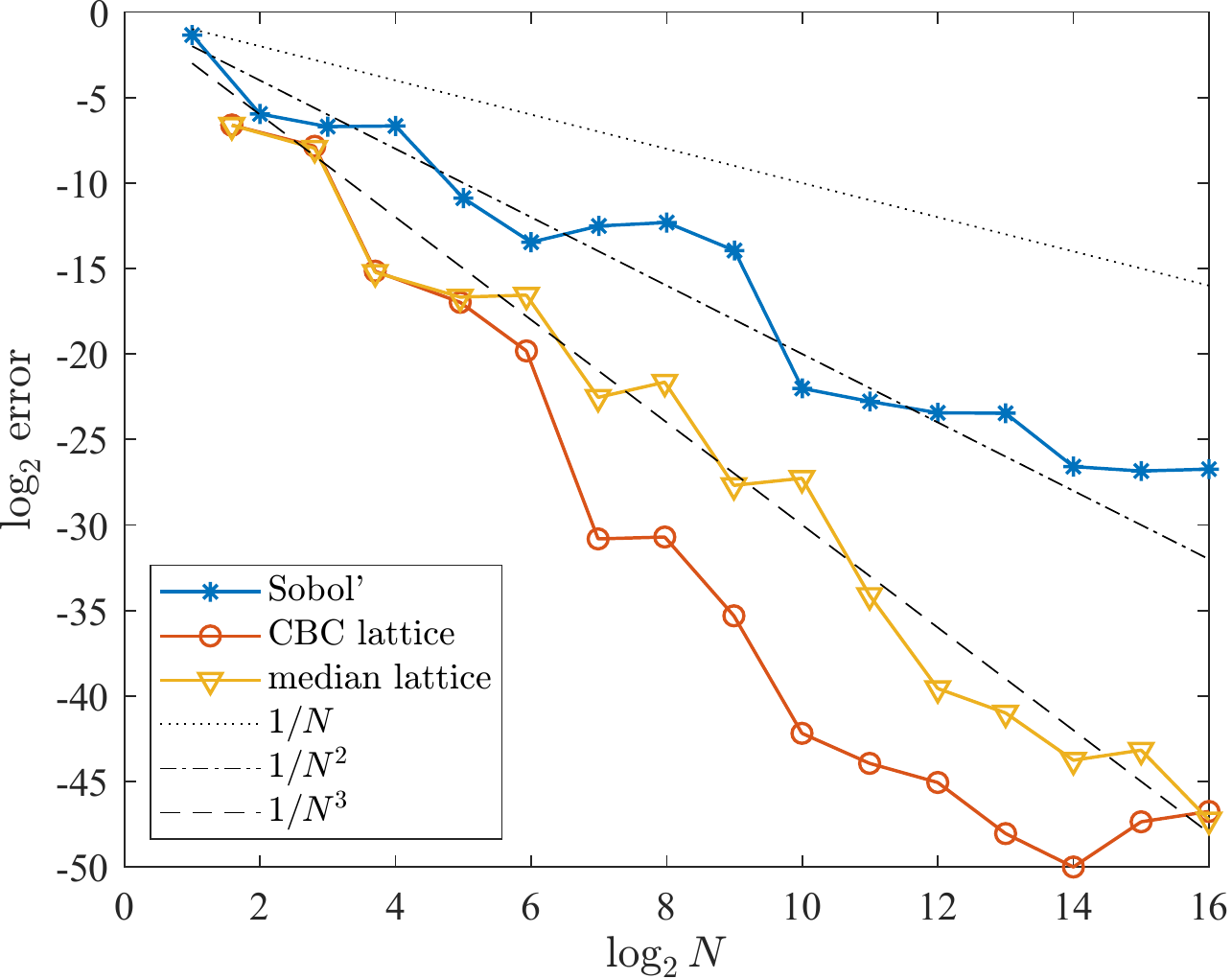}
    \includegraphics[width=0.45\linewidth]{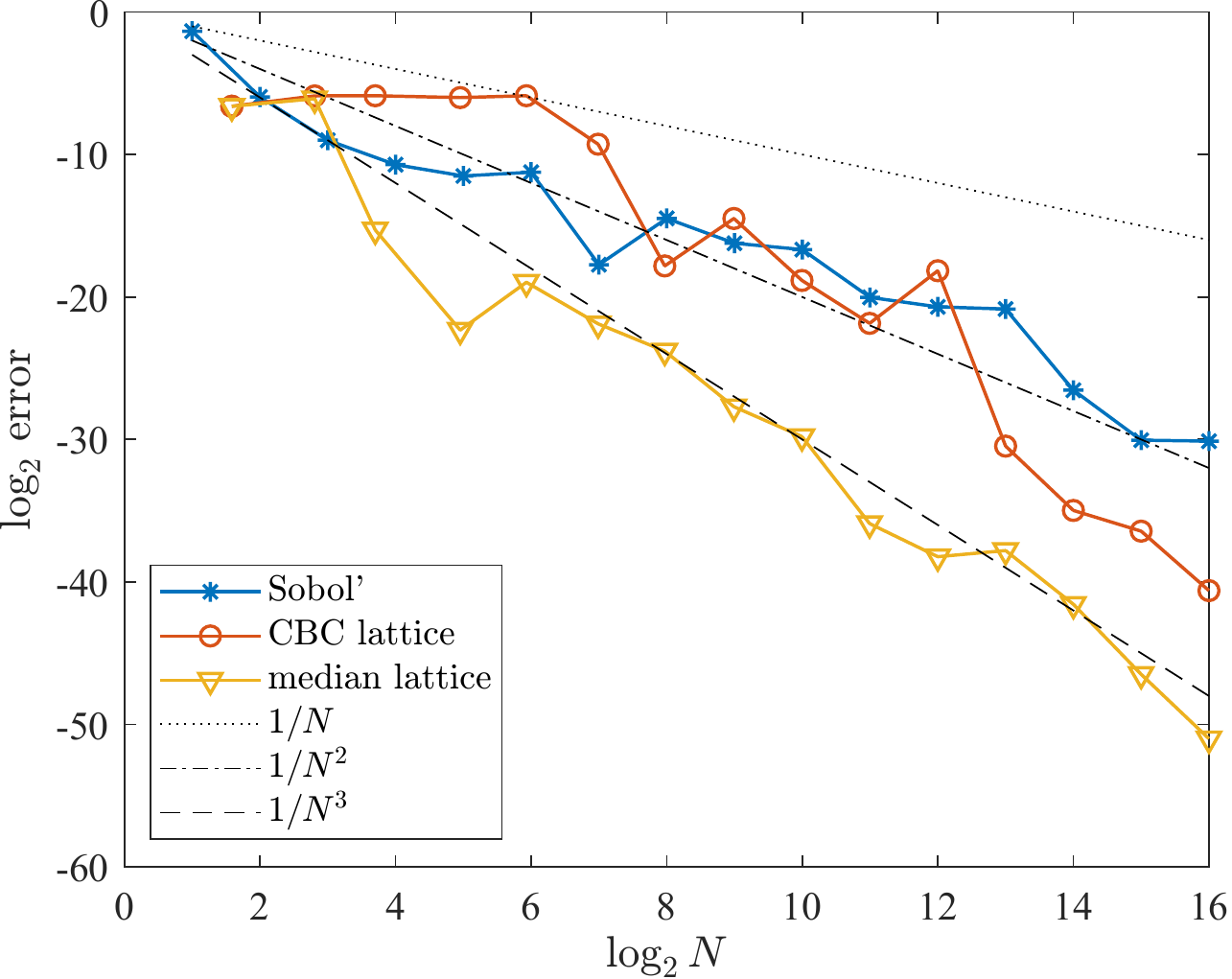}\\
    \caption{Comparison of the integration error by our median lattice rule (yellow), QMC rule using Sobol' points (blue), and rank-1 lattice rule with the fast CBC algorithm (orange). 
	The results are shown for the test function $f^{\per}_{\beta,\bsomega}$ with the choices $\beta=2$ and $\omega_j=1/j^3$ (upper left), $\beta=2$ and $\omega_j=1/(s-j+1)^3$ (upper right), $\beta=5$ and $\omega_j=1/j^6$ (lower left), and $\beta=5$ and $\omega_j=1/(s-j+1)^6$ (lower right).}
    \label{fig:lattice}
\end{figure}

To show that our median lattice rule performs well for functions of non-product forms, let us consider the additional test functions given by
\[ f^{\per,\mathrm{cyc}}_{\beta}(\bsx)=\frac{1}{5}\sum_{\ell=1}^{5}\prod_{j=1}^{s/5}g_{\beta}(x_{j+s(\ell-1)/5})\quad
\text{and}\quad f^{\per,\mathrm{mod}}_{\beta}(\bsx)=\frac{1}{5}\sum_{\ell=1}^{5}\prod_{j=1}^{s/5}g_{\beta}(x_{\ell+5(j-1)}),\]
respectively, with $\beta=5$ and $s=20$. We have that $I_s(f^{\per,\mathrm{cyc}}_{\beta})=I_s(f^{\per,\mathrm{mod}}_{\beta})=1$ and these two integrands belong to the Korobov space with $\alpha=\beta$. The results for the two integrands are shown in the corresponding panels of Figure~\ref{fig:lattice_nonproduct}. For large $N$, both our median lattice rule and the rank-1 lattice rule constructed by the CBC algorithm are superior to the QMC rule using the Sobol' points. Although the difference between $f^{\per,\mathrm{cyc}}_{\beta}$ and $f^{\per,\mathrm{mod}}_{\beta}$ lies only in the ordering of variables, the convergence behavior of the rank-1 lattice rule constructed by the CBC algorithm is not consistent for these functions and a strange zig-zag pattern shows up for $f^{\per,\mathrm{mod}}_{\beta}$. On the contrary, our median lattice rule is not subject to the difference between the ordering of variables and performs almost equivalently.

\begin{figure}
    \centering
    \includegraphics[width=0.45\linewidth]{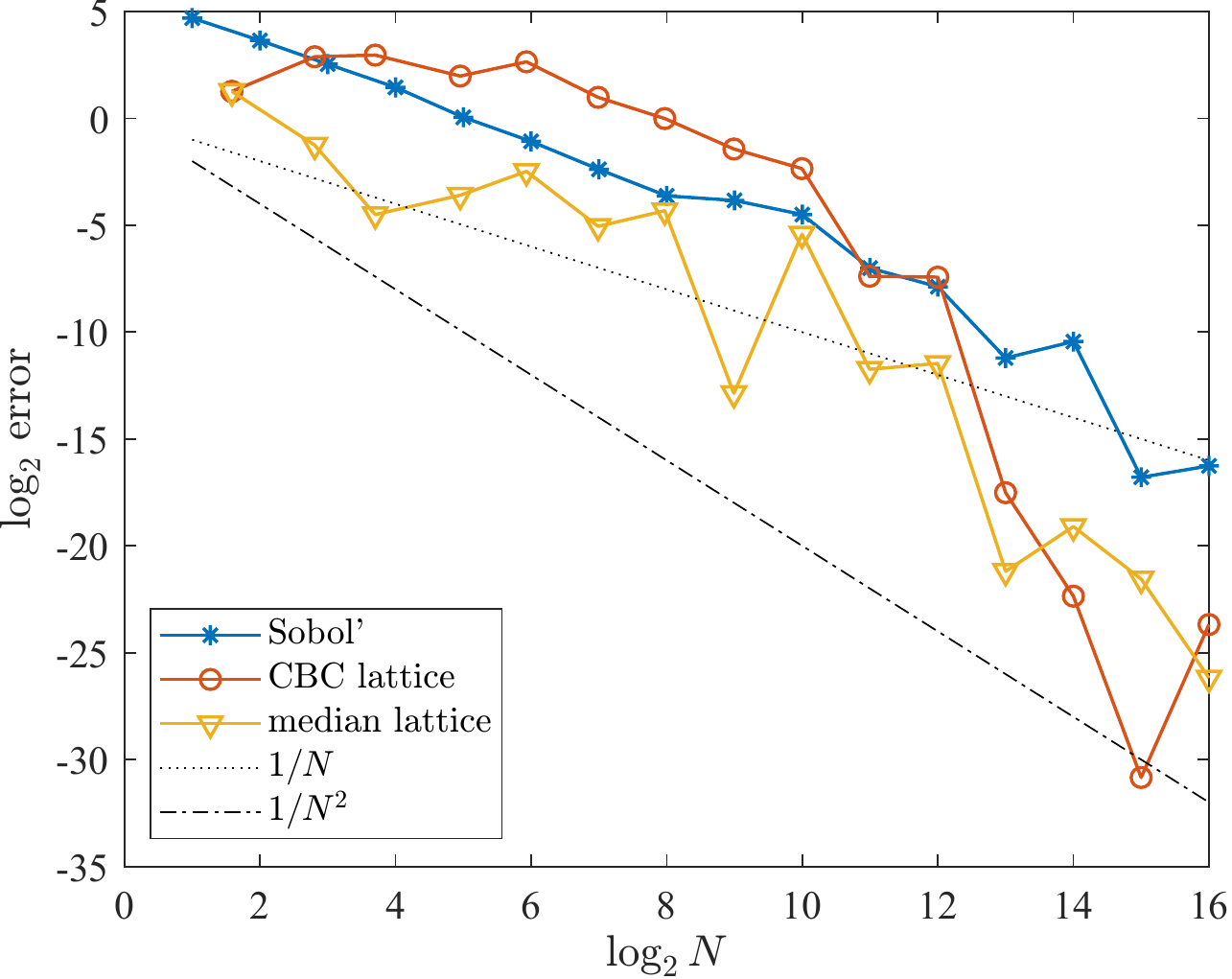}
    \includegraphics[width=0.45\linewidth]{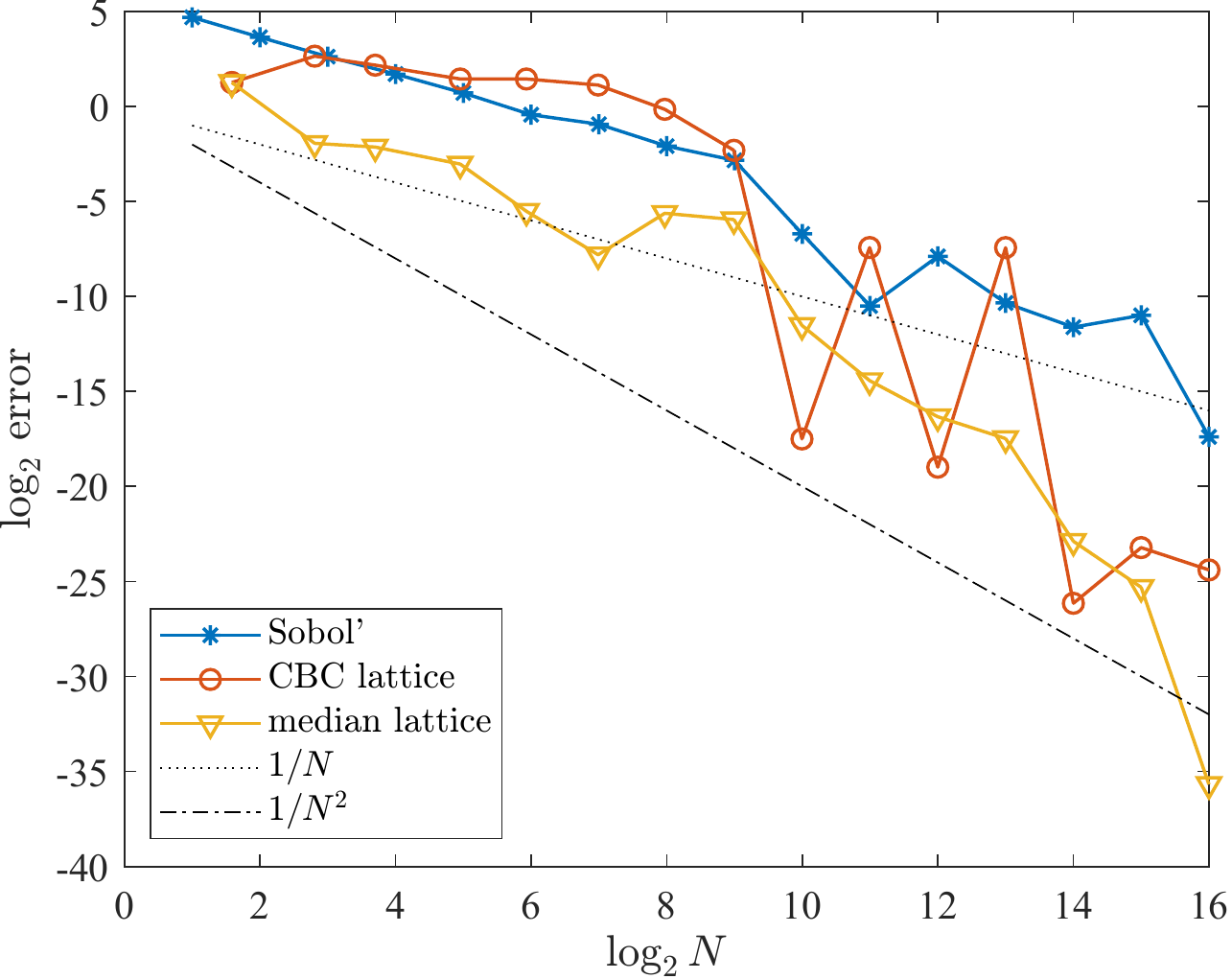}\\
    \caption{Comparison of the integration error by our median lattice rule (yellow), QMC rule using Sobol' points (blue), and rank-1 lattice rule with the fast CBC algorithm (orange). 
	The results are shown for the test functions $f^{\per,\mathrm{cyc}}_{\beta}$ (left) and $f^{\per,\mathrm{mod}}_{\beta}$ (right) with the choices $\beta=5$ and $s=20$.}
    \label{fig:lattice_nonproduct}
\end{figure}

\end{example}

\subsection{High-order polynomial lattice rules for non-periodic functions}

Our next examples concern high-order polynomial lattice rules. Here we fix the precision to $n=52$ and always use the primitive polynomial $p(x)=x^{52}+x^3+1$, found in \cite{HM92}, as the modulus of the polynomial lattice point sets.

\begin{figure}
    \centering
    \includegraphics[width=0.45\linewidth]{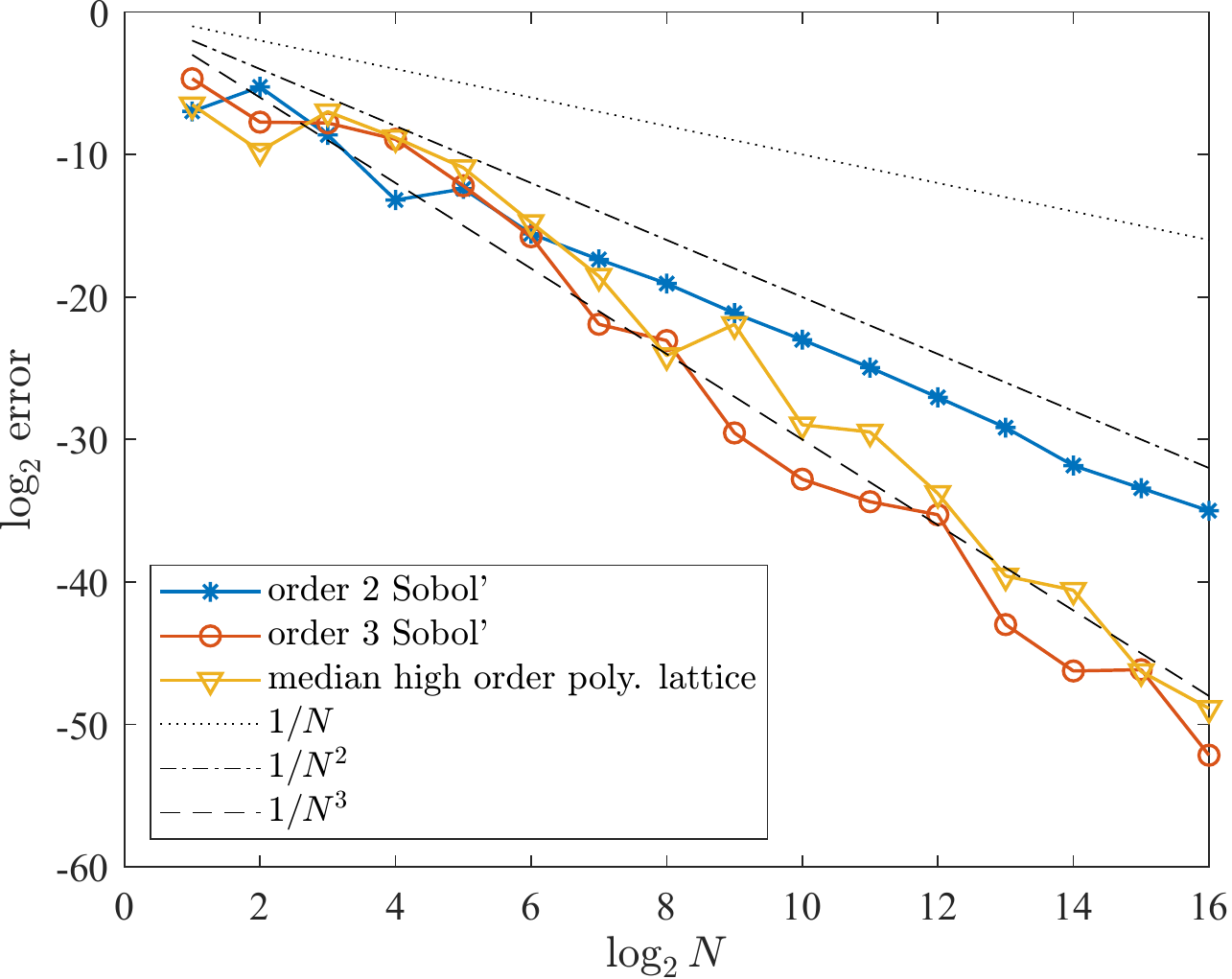}
    \includegraphics[width=0.45\linewidth]{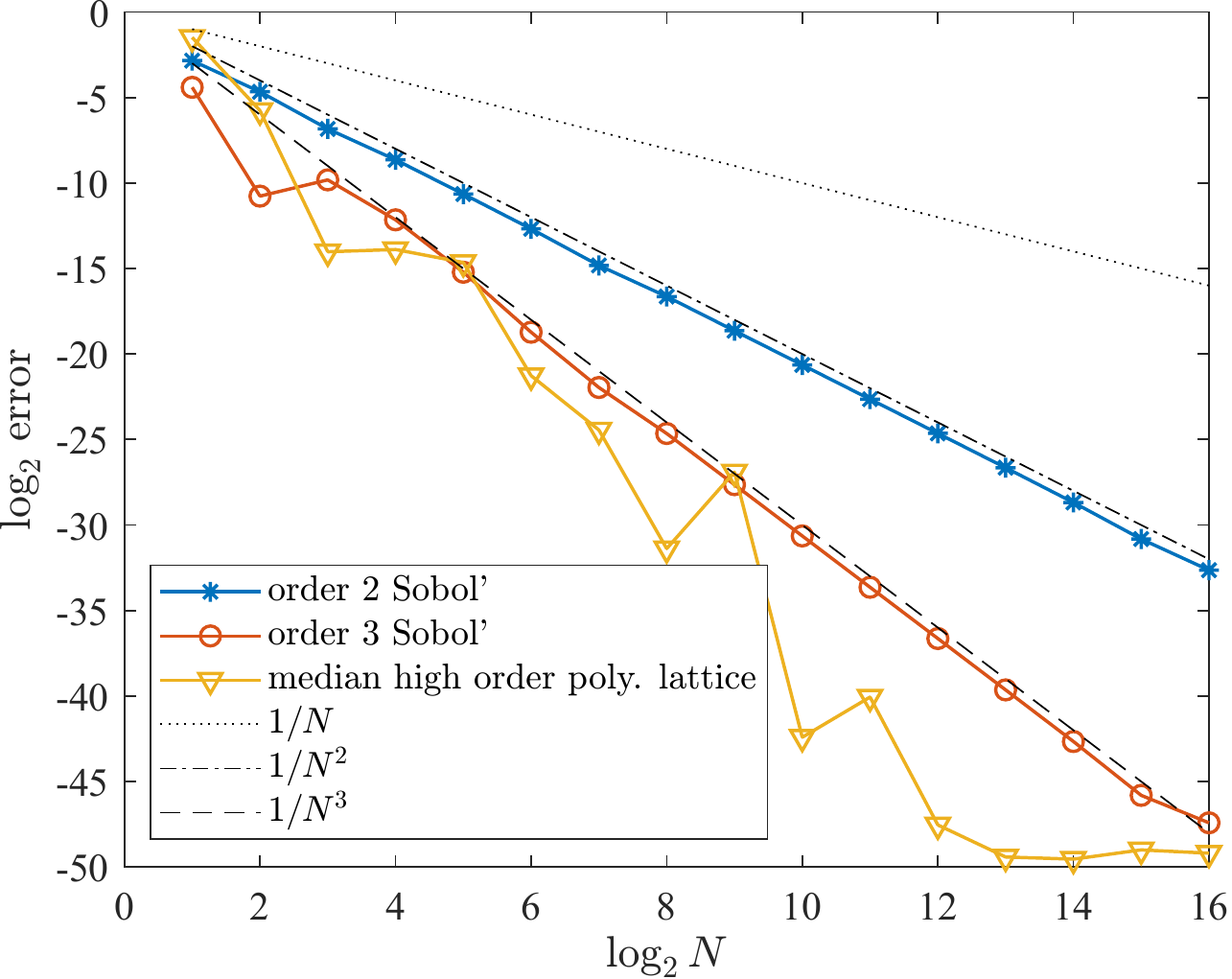}
    \caption{Comparison of the one-dimensional integration error by our median high-order polynomial lattice rule (yellow), QMC rules using order 2 Sobol' points (blue) and order 3 Sobol' points (orange). The results are shown for the test functions $f^{\nonper}_{1}$ (left) and $f^{\nonper}_{2}$ (right), respectively.}
    \label{fig:HOPL_1d}
\end{figure}

\begin{example}\rm 
We first consider the two following one-dimensional test functions:
\[ 
  f^{\nonper}_{1}(x)=x^3 (1/4+\log x)\qquad \text{and}\qquad  f^{\nonper}_{2}(x)=xe^{x/4}. 
\]
We can see that the third derivative of $f^{\nonper}_{1}$ is in $L_q([0,1))$ for any $1\leq q< \infty$, whereas the fourth derivative is not in $L_1([0,1))$, implying that $f^{\nonper}_{1}\in \cF^{\Sob}_{1,3,\gamma,q}$ but $f^{\nonper}_{1}\not\in \cF^{\Sob}_{1,4,\gamma,1}$. Thus $f^{\nonper}_{1}$ has a finite smoothness. On the other hand, $f^{\nonper}_{2}$ is obviously infinitely differentiable, so that $f^{\nonper}_{2}\in \cF^{\Sob}_{1,\alpha,\gamma,q}$ for any $\alpha\geq 2$ and $1\leq q\leq \infty$. Note that $I_1(f^{\nonper}_{1})=0$ and $I_1(f^{\nonper}_{2})=16-12e^{1/4}$. We compare our median high-order polynomial lattice rule with $r=11$ and QMC rules using order 2 and order 3 Sobol' points
constructed by the interlacing procedure of \cite{Dick08}, 
with the direction numbers provided in MATLAB (taken from \cite{JK03}). 
To construct a Sobol' point set of order $d$ by interlacing, we first construct a $ds$-dimensional Sobol point set 
with $2^m$ points (with $s=1$ in this case) and then apply the digit interlacing procedure defined in \cite{Dick08}
to obtain the digits of the $s$-dimensional points.  This procedure extracts the first $m$ digits of the 
$ds$-dimensional points and reorders them in a special way to obtain the first $dm$ digits of the $s$-dimensional points.

The results for the two one-dimensional functions are shown in Figure~\ref{fig:HOPL_1d}, respectively. As we can see from the result for $f^{\nonper}_{1}$, the QMC rule using order 2 Sobol' points cannot fully exploit the smoothness of the function and the error decays at the rate of $N^{-2}$. On the other hand, the QMC rule using order 3 Sobol' points and our median high-order polynomial lattice rule can exploit the smoothness and achieves the convergence rate of $N^{-3}$. For the infinitely differentiable function $f^{\nonper}_{2}$, the plot suggests that our median high-order polynomial lattice rule may converge even faster than $N^{-3}$. 
These numerical results show the major advantage of our proposed rule in terms of adaptivity in smoothness.
\end{example}

\begin{figure}
    \centering
    \includegraphics[width=0.45\linewidth]{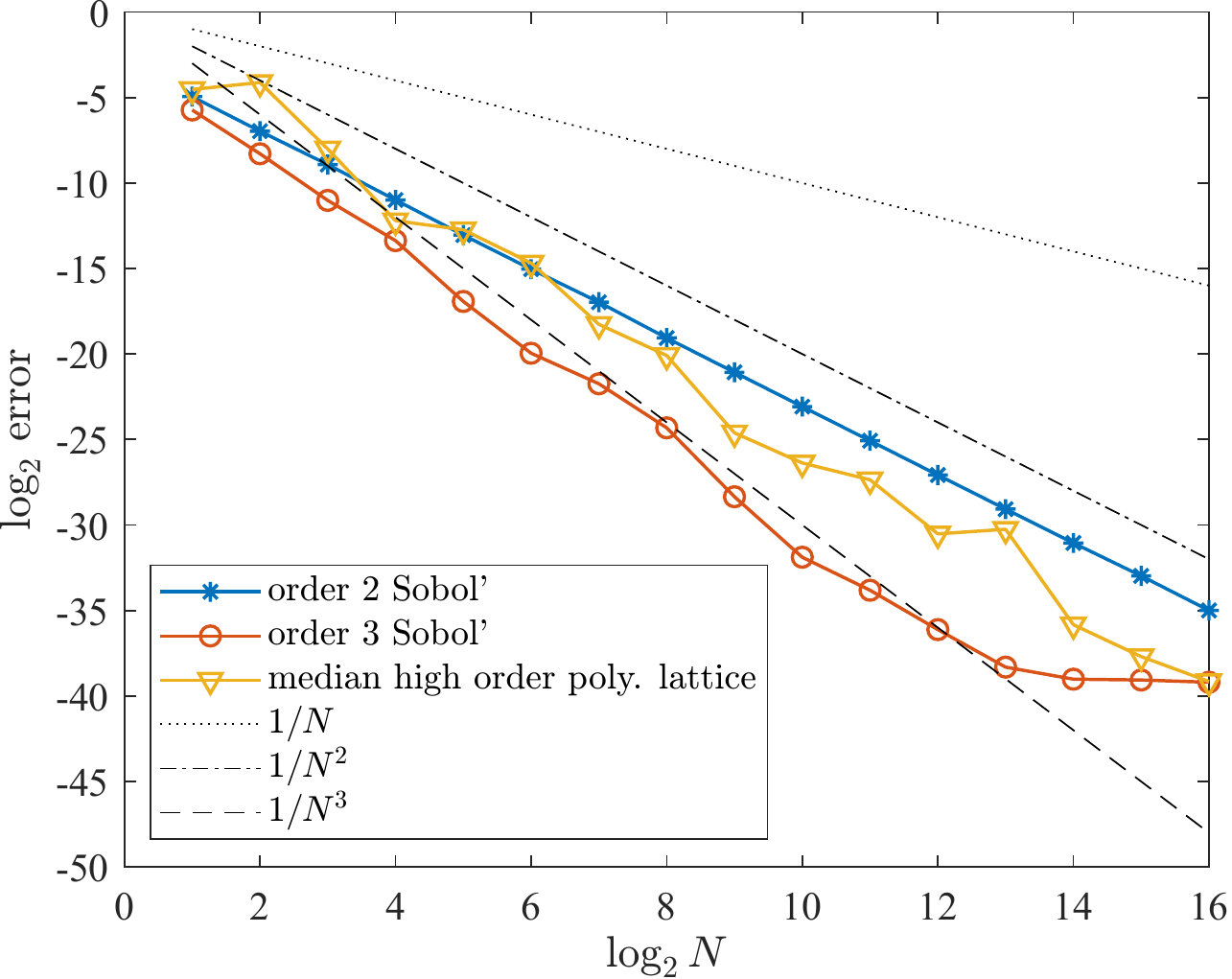}
    \includegraphics[width=0.45\linewidth]{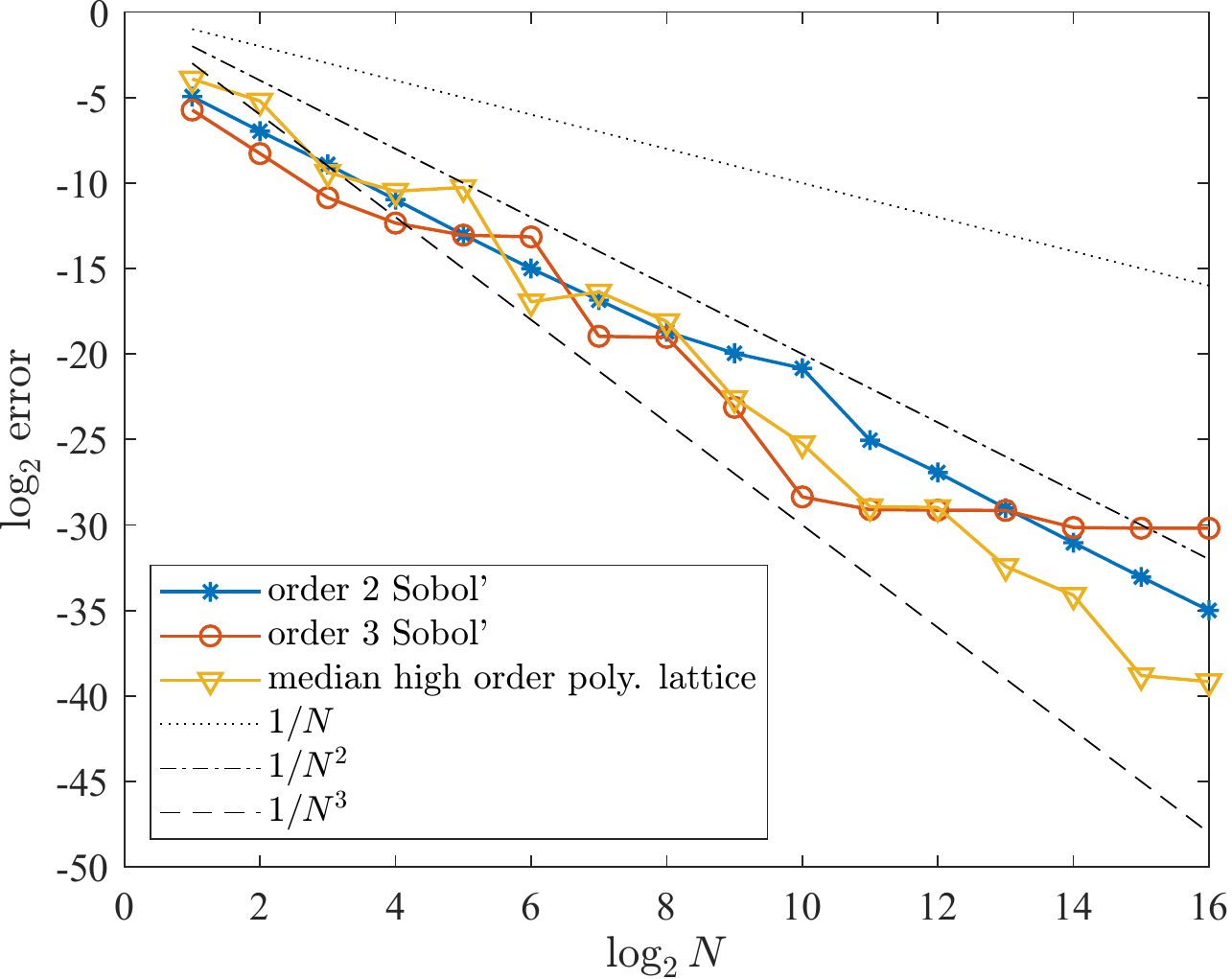}
    \caption{Comparison of the multi-dimensional integration error by our median high-order polynomial lattice rule (yellow), QMC rules using order 2 Sobol' points (blue) and order 3 Sobol' points (orange). The results are shown for the test functions $f^{\nonper}_{3,\bsomega}$ (left) and $f^{\nonper}_{3,\bsomega,\mathrm{flip}}$ (right) with the choice $\omega_j=1/(4j^4)$}.
    \label{fig:HOPL_10d}
\end{figure}

\begin{example}\rm 
Finally, we consider the two multivariate non-periodic test functions
\[ 
  f^{\nonper}_{3,\bsomega}(\bsx) = \exp\left( -\sum_{j=1}^{s}\omega_jx_j\right)
	 \qquad \text{and}\qquad 
	f^{\nonper}_{3,\bsomega,\mathrm{flip}}(\bsx)=\exp\left( -\sum_{j=1}^{s}\omega_jx_{s-j+1}\right),
\]
with $s=10$ and $\omega_j=1/(4j^4)$. It is obvious that both $f^{\nonper}_{3,\bsomega}$ and $f^{\nonper}_{3,\bsomega,\mathrm{flip}}$ are infinitely differentiable and belong to $\cF^{\Sob}_{s,\alpha,\bsgamma,q}$ with arbitrary $\alpha\geq 2$ and $1\leq q\leq \infty$. Note that $f^{\nonper}_{3,\bsomega,\mathrm{flip}}$ is defined by reordering the variables of $f^{\nonper}_{3,\bsomega}$ so that $x_j$ is replaced by $x_{s-j+1}$, and that we have
\[ 
  I_s(f^{\nonper}_{3,\bsomega})=I_s(f^{\nonper}_{3,\bsomega,\mathrm{flip}})=\prod_{j=1}^{s}\frac{1-\exp(-\omega_j)}{\omega_j}.
\]
The variables are ordered by decreasing order of importance in the first function, and by increasing order in the second one. As our median high-order polynomial lattice rule, based on random choices of generating vectors, does not care about the ordering of variables, it should perform the same for $f^{\nonper}_{3,\bsomega}$ and $f^{\nonper}_{3,\bsomega,\mathrm{flip}}$.

The results are shown in Figure~\ref{fig:HOPL_10d}. Here again, we compare our median high-order polynomial lattice rule with $r=11$ and QMC rules using order 2 and order 3 interlaced Sobol' points. 
For the function $f^{\nonper}_{3,\bsomega}$, our median high-order polynomial lattice rule can exploit the smoothness better than the QMC rule using order 2 Sobol' points. 
The QMC rule using order 3 Sobol' points exploits the smoothness of the integrand best and the error decays at the rate of $N^{-3}$ and outperforms our median high-order polynomial lattice rule approximately by a constant factor for small $N$, 
but this rate breaks down at around $\log_2 N = 13$ and our median rule catches up at $\log_2 N = 16$.

For $f^{\nonper}_{3,\bsomega,\mathrm{flip}}$, the situation changes. Our median high-order polynomial lattice rule is now comparable to the QMC rule using order 3 Sobol' points when $N$ is small, and performs better for larger $N$. The error decays approximately at the rate of $N^{-2.5}$, which the QMC rule using order 2 Sobol' points cannot attain. The slowdown of the convergence for the QMC rule using order 3 Sobol' points might be due to the misspecification of important variables. Although $x_s,x_{s-1},\ldots$ are the order of the relatively important variables for $f^{\nonper}_{3,\bsomega,\mathrm{flip}}$, we use the later coordinates of order 3 Sobol' points, whose lower-dimensional projections are not well-distributed compared to the earlier coordinates. In this sense, the median high-order polynomial lattice rule is more robust and adaptive to the integrand at hand.
\end{example}

\section*{Acknowledgments}

The authors would like to thank Art Owen and Zexin Pan for sharing their preprint and having useful discussions, and also thank Mario Ullrich for giving useful comments.

\bibliographystyle{plain}
\bibliography{median-lattice-ref.bib}

\begin{thebibliography}{10}

\bibitem{BDLNP12}
J.~Baldeaux, J.~Dick, G.~Leobacher, D.~Nuyens, and F.~Pillichshammer.
\newblock Efficient calculation of the worst-case error and (fast)
  component-by-component construction of higher order polynomial lattice rules.
\newblock {\em Numerical Algorithms}, 59(3):403--431, 2012.

\bibitem{CKNS16}
R.~Cools, F.~Y. Kuo, D.~Nuyens, and G.~Suryanarayana.
\newblock Tent-transformed lattice rules for integration and approximation of
  multivariate non-periodic functions.
\newblock {\em Journal of Complexity}, 36:166--181, 2016.

\bibitem{Dick08}
J.~Dick.
\newblock Walsh spaces containing smooth functions and quasi–{M}onte {C}arlo
  rules of arbitrary high order.
\newblock {\em SIAM Journal on Numerical Analysis}, 46(3):1519--1553, 2008.

\bibitem{DG21}
J.~Dick and T.~Goda.
\newblock Stability of lattice rules and polynomial lattice rules constructed
  by the component-by-component algorithm.
\newblock {\em Journal of Computational and Applied Mathematics}, 382:113062,
  2021.

\bibitem{DHP15}
J.~Dick, A.~Hinrichs, and F.~Pillichshammer.
\newblock Proof techniques in quasi-{M}onte {C}arlo theory.
\newblock {\em Journal of Complexity}, 31:327--371, 2015.

\bibitem{DKLNS14}
J.~Dick, F.~Y. Kuo, Q.~T.~Le Gia, D.~Nuyens, and C.~Schwab.
\newblock Higher order {Q}{M}{C} {P}etrov-{G}alerkin discretization for affine
  parametric operator equations with random field inputs.
\newblock {\em SIAM Journal on Numerical Analysis}, 52(6):2676--2702, 2014.

\bibitem{DKS13}
J.~Dick, F.~Y. Kuo, and I.~H. Sloan.
\newblock High-dimensional integration: the quasi-{M}onte {C}arlo way.
\newblock {\em Acta Numerica}, 22:133--288, 2013.

\bibitem{DNP14}
J.~Dick, D.~Nuyens, and F.~Pillichshammer.
\newblock Lattice rules for nonperiodic smooth integrands.
\newblock {\em Numerische Mathematik}, 126:259--291, 2014.

\bibitem{DP07}
J.~Dick and F.~Pillichshammer.
\newblock Strong tractability of multivariate integration of arbitrary high
  order using digitally shifted polynomial lattice rules.
\newblock {\em Journal of Complexity}, 23:436--453, 2007.

\bibitem{DickPilli10}
J.~Dick and F.~Pillichshammer.
\newblock {\em Digital Nets and Sequences: Discrepancy Theory and Quasi-Monte
  Carlo Integration}.
\newblock Cambridge University Press, Cambridge, 2010.

\bibitem{DSWW06}
J.~Dick, I.~H. Sloan, X.~Wang, and H.~Wo\'{z}niakowski.
\newblock Good lattice rules in weighted {K}orobov spaces with general weights.
\newblock {\em Numerische Mathematik}, 103(1):63--97, 2006.

\bibitem{EKNO21}
A.~Ebert, P.~Kritzer, D.~Nuyens, and O.~Osisiogu.
\newblock Digit-by-digit and component-by-component constructions of lattice
  rules for periodic functions with unknown smoothness.
\newblock {\em Journal of Complexity}, 66:101555, 2021.

\bibitem{G15}
T.~Goda.
\newblock Good interlaced polynomial lattice rules for numerical integration in
  weighted {W}alsh spaces.
\newblock {\em Journal of Computational and Applied Mathematics}, 285:279--294,
  2015.

\bibitem{G16}
T.~Goda.
\newblock Quasi-{M}onte {C}arlo integration using digital nets with
  antithetics.
\newblock {\em Journal of Computational and Applied Mathematics}, 304:26--42,
  2016.

\bibitem{GD15}
T.~Goda and J.~Dick.
\newblock Construction of interlaced scrambled polynomial lattice rules of
  arbitrary high order.
\newblock {\em Foundations of Computational Mathematics}, 15(5):1245--1278,
  2015.

\bibitem{GSY19}
T.~Goda, K.~Suzuki, and T.~Yoshiki.
\newblock Lattice rules in non-periodic subspaces of {S}obolev spaces.
\newblock {\em Numerische Mathematik}, 141:399--427, 2019.

\bibitem{HM92}
T.~Hansen and G.~L. Mullen.
\newblock Primitive polynomials over finite fields.
\newblock {\em Mathematics of Computation}, 59:639--643, 1992.

\bibitem{Hickernell02}
F.~J. Hickernell.
\newblock Obtaining ${O}(n^{-2+\epsilon})$ convergence for lattice quadrature
  rules.
\newblock In K.-T. Fang, H.~Niederreiter, and F.~J. Hickernell, editors, {\em
  Monte Carlo and Quasi-Monte Carlo Methods 2000}, pages 274--289, Berlin,
  2002. Springer-Verlag.

\bibitem{JK03}
S.~Joe and F.~Y. Kuo.
\newblock Remark on algorithm 659: {I}mplementing {S}obol's quasirandom
  sequence generator.
\newblock {\em ACM Transactions on Mathematical Software}, 29(1):49--57, 2003.

\bibitem{Korobov59}
N.~M. Korobov.
\newblock The approximate computation of multiple integrals (in {R}ussian).
\newblock {\em Doklady Akademii Nauk SSSR}, 124:1207--–1210, 1959.

\bibitem{KKNU19}
P.~Kritzer, F.~Y. Kuo, D.~Nuyens, and M.~Ullrich.
\newblock Lattice rules with random $n$ achieve nearly the optimal
  $\mathcal{O}(n^{-\alpha-1/2})$ error independently of the dimension.
\newblock {\em Journal of Approximation Theory}, 240:96--113, 2019.

\bibitem{Kuo03}
F.~Y. Kuo.
\newblock Component-by-component constructions achieve the optimal rate of
  convergence for multivariate integration in weighted {K}orobov and {S}obolev
  spaces.
\newblock {\em Journal of Complexity}, 19:301--320, 2003.

\bibitem{KJ02}
F.~Y. Kuo and S.~Joe.
\newblock Component-by-component construction of good lattice rules with a
  composite number of points.
\newblock {\em Journal of Complexity}, 18:943--976, 2002.

\bibitem{KSS12}
F.~Y. Kuo, C.~Schwab, and I.~H. Sloan.
\newblock Quasi-{M}onte {C}arlo finite element methods for a class of elliptic
  partial differential equations with random coefficients.
\newblock {\em SIAM Journal on Numerical Analysis}, 50:3351--3374, 2012.

\bibitem{vLEC18a}
P.~L'Ecuyer.
\newblock Randomized quasi-{Monte Carlo}: An introduction for practitioners.
\newblock In P.~W. Glynn and A.~B. Owen, editors, {\em {M}onte {C}arlo and
  Quasi-{M}onte {C}arlo Methods: {MCQMC 2016}}, pages 29--52, Berlin, 2018.
  Springer.

\bibitem{vLEC02a}
P.~L'Ecuyer and C.~Lemieux.
\newblock Recent advances in randomized quasi-{M}onte {C}arlo methods.
\newblock In M.~Dror, P.~L'Ecuyer, and F.~Szidarovszky, editors, {\em Modeling
  Uncertainty: An Examination of Stochastic Theory, Methods, and Applications},
  pages 419--474. Kluwer Academic, Boston, 2002.

\bibitem{LatNetbuilder}
P.~L'Ecuyer, P.~Marion, M.~Godin, and F.~Puchhammer.
\newblock A tool for custom construction of {Q}{M}{C} and {R}{Q}{M}{C} point
  sets, December 2020.
\newblock arXiv:2012.10263.

\bibitem{vLEC12a}
P.~L'Ecuyer and D.~Munger.
\newblock On figures of merit for randomly-shifted lattice rules.
\newblock In L.~Plaskota and H.~Wo\'{z}niakowski, editors, {\em {M}onte {C}arlo
  and Quasi-{M}onte {C}arlo Methods 2010}, pages 133--159, Berlin, 2012.
  Springer-Verlag.

\bibitem{Latticebuilder}
P.~L'Ecuyer and D.~Munger.
\newblock Algorithm 958: {L}attice {B}uilder: {A} general software tool for
  constructing rank-1 lattice rules.
\newblock {\em ACM Transactions on Mathematical Software}, 42(2):15, 2016.

\bibitem{Lemieux09}
C.~Lemieux.
\newblock {\em Monte Carlo and Quasi-Monte Carlo Sampling}.
\newblock Springer, New York, 2009.

\bibitem{vLEM03a}
C.~Lemieux and P.~L'Ecuyer.
\newblock Randomized polynomial lattice rules for multivariate integration and
  simulation.
\newblock {\em {SIAM} Journal on Scientific Computing}, 24(5):1768--1789, 2003.

\bibitem{LeobPilli14}
G.~Leobacher and F.~Pillichshammer.
\newblock {\em Introduction to Quasi-Monte Carlo Integration and Applications}.
\newblock Birkh\"{a}user, Cham, 2014.

\bibitem{Merkle05}
M.~Merkle.
\newblock Jensen’s inequality for medians.
\newblock {\em Statistics \& Probability Letters}, 71(3):277--281, 2005.

\bibitem{vNIC14a}
J.~A. Nichols and F.~Y. Kuo.
\newblock Fast {CBC} construction of randomly shifted lattice rules achieving
  ${O}(n^{-1+\delta})$ convergence rate for unbounded integrands over
  $\mathbb{R}^s$ in weighted spaces with {POD} weights.
\newblock {\em Journal of Complexity}, 30:444--468, 2014.

\bibitem{Niederreiter92b}
H.~Niederreiter.
\newblock Low-discrepancy point sets obtained by digital constructions over
  finite fields.
\newblock {\em Czechoslovak Mathematical Journal}, 42:143--166, 1992.

\bibitem{Niederreiter92}
H.~Niederreiter.
\newblock {\em Random Number Generation and Quasi-Monte Carlo Methods}.
\newblock SIAM, Philadelphia, 1992.

\bibitem{NW_Tract1}
E.~Novak and H.~Wo\'{z}niakowski.
\newblock {\em Tractability of Multivariate Problems. Volume I: Linear
  Information}.
\newblock EMS, Z\"{u}rich, 2008.

\bibitem{iNUY20m}
D.~Nuyens.
\newblock The magic point shop, 2020.
\newblock \url{https://people.cs.kuleuven.be/~dirk.nuyens/qmc-generators/}.

\bibitem{NC06}
D.~Nuyens and R.~Cools.
\newblock Fast algorithms for component-by-component construction of rank-1
  lattice rules in shift-invariant reproducing kernel {H}ilbert spaces.
\newblock {\em Mathematics of Computation}, 75:903--920, 2006.

\bibitem{PanOwen2021}
Z.~Pan and A.~B. Owen.
\newblock Super-polynomial accuracy of one dimensional randomized nets using
  the median-of-means, November 2021.
\newblock arXiv:2111.12676.

\bibitem{Pillichshammer12}
F.~Pillichshammer.
\newblock Polynomial lattice point sets.
\newblock In L.~Plaskota and H.~Wo\'{z}niakowski, editors, {\em Monte Carlo and
  Quasi-Monte Carlo Methods 2010}, pages 189--210, Berlin, 2012.
  Springer-Verlag.

\bibitem{RS62}
J.~B. Rosser and L.~Schoenfeld.
\newblock Approximate formulas for some functions of prime numbers.
\newblock {\em Illinois Journal of Mathematics}, 6(1):64--94, 1962.

\bibitem{SloanJoe94}
I.~H. Sloan and S.~Joe.
\newblock {\em Lattice Methods for Multiple Integration}.
\newblock Oxford Science Publications, New York, 1994.

\bibitem{SR02}
I.~H. Sloan and A.~V. Reztsov.
\newblock Component-by-component construction of good lattice rules.
\newblock {\em Mathematics of Computation}, 71:263--273, 2002.

\bibitem{SW98}
I.~H. Sloan and H.~Wo\'{z}niakowski.
\newblock When are quasi-{M}onte {C}arlo algorithms efficient for
  high-dimensional integrals?
\newblock {\em Journal of Complexity}, 14:1--33, 1998.

\end{thebibliography}

\end{document}